\documentclass[a4paper,12pt]{amsart}

\usepackage[a4paper, top=3cm, left=2.5cm, right=2.5cm, bottom=3cm]{geometry}
\usepackage{amsmath, amsthm, pb-diagram}
\usepackage{amssymb}
\usepackage{amsfonts,graphicx,color}
\usepackage{enumerate, fancyhdr, dsfont}
\usepackage[normalem]{ulem}

    \newcommand{\Cwf}{\mathcal{C}}

    \newcommand{\Hwf}{\mathcal{H}}
    \newcommand{\Iwf}{\mathcal{I}}
    
    \newcommand{\Mwf}{\mathcal{M}}
    \newcommand{\Nwf}{\mathcal{N}}

    \newcommand{\Swf}{\mathcal{S}}

    \newcommand{\afrak}{\mathfrak{a}}
    \newcommand{\bfrak}{\mathfrak{b}}
    \newcommand{\cfrak}{\mathfrak{c}}
    \newcommand{\dfrak}{\mathfrak{d}}

    \newcommand{\pfrak}{\mathfrak{p}}
    \newcommand{\rfrak}{\mathfrak{r}}
    \newcommand{\sfrak}{\mathfrak{s}}
    
    \newcommand{\ufrak}{\mathfrak{u}}

    \newcommand{\menos}{\smallsetminus}

    \newcommand{\frestr}{\!\!\upharpoonright\!\!}

    \newcommand{\add}{\mbox{\rm add}}
    \newcommand{\cov}{\mbox{\rm cov}}
    \newcommand{\non}{\mbox{\rm non}}
    \newcommand{\cof}{\mbox{\rm cof}}
    \newcommand{\limdir}{\mbox{\rm limdir}}
    
    \newcommand{\ral}{\mathfrak{R}}
     \newcommand{\M}{\mathcal{M}}
    
    \newcommand{\Bor}{\mathds{B}}
    \newcommand{\Cor}{\mathds{C}}
    \newcommand{\Dor}{\mathds{D}}
    \newcommand{\Eor}{\mathds{E}}
    \newcommand{\Hor}{\mathds{H}}
    
    \newcommand{\Loc}{\mathds{LOC}}
    \newcommand{\Mor}{\mathds{M}}
    \newcommand{\Por}{\mathds{P}}
    \newcommand{\Qor}{\mathds{Q}}
    
    \newcommand{\Sor}{\mathds{S}}
    
    \newcommand{\Qnm}{\dot{\mathds{Q}}}

    \newcommand{\R}{\mathbb{R}}

    \newcommand{\cf}{\mbox{\rm cf}}

    \newcommand{\wdt}{\mathrm{width}}

    \newcommand{\la}{\langle}
    \newcommand{\ra}{\rangle}
    \newcommand{\PP}{\mathds{P}}
    \newcommand{\QQ}{\mathds{Q}}

   \newcommand{\Lim}{\mathrm{Lim}}

\title[Coherent systems of finite support iterations]{Coherent systems of finite support iterations}

\author{Vera Fischer}
\address{Kurt G\"odel Research Center, University of Vienna, W\"ahringer Strasse 25, 1090 Vienna, Austria}
\email{vera.fischer@univie.ac.at}
\urladdr{http://www.logic.univie.ac.at/~vfischer/}

\author{Sy D. Friedman}
\address{Kurt G\"odel Research Center, University of Vienna, W\"ahringer Strasse 25, 1090 Vienna, Austria}
\email{sdf@univie.ac.at}
\urladdr{http://www.logic.univie.ac.at/~sdf/}

\author{Diego A. Mej\'ia}
\address{Creative Science Course (Mathematics), Faculty of Science, Shizuoka University, Ohya 836, Suruga-ku, Shizuoka-shi, Japan 422-8529.}
\email{diego.mejia@shizuoka.ac.jp}
\urladdr{http://www.researchgate.com/profile/Diego\_Mejia2}

\author{Diana C. Montoya}
\address{Kurt G\"odel Research Center, University of Vienna, W\"ahringer Strasse 25, 1090 Vienna, Austria}
\email{dcmontoyaa@gmail.com}
\urladdr{http://www.logic.univie.ac.at/~montoyd8/}

\thanks{The first author would like to thank the Austrian Science Fund (FWF) for the generous support through Lise-Meitner Project M1365-N13. The second and fourth authors were supported by the Austrian
  Science Fund (FWF) Project P25748. The third author was supported by the
  FWF Projects P23875-N13 and I1272-N25}



\begin{document}

\makeatletter
\def\@roman#1{\romannumeral #1}
\makeatother

\theoremstyle{plain}
  \newtheorem{theorem}{Theorem}[section]
  \newtheorem{corollary}[theorem]{Corollary}
  \newtheorem{lemma}[theorem]{Lemma}
  \newtheorem{prop}[theorem]{Proposition}
  \newtheorem{clm}[theorem]{Claim}
  \newtheorem{exer}[theorem]{Exercise}
  \newtheorem{question}[theorem]{Question}
  \newtheorem{problem}[theorem]{Problem}
  \newtheorem*{thm}{Theorem}
\theoremstyle{definition}
  \newtheorem{definition}[theorem]{Definition}
  \newtheorem{example}[theorem]{Example}
  \newtheorem{remark}[theorem]{Remark}
  \newtheorem{notation}[theorem]{Notation}
  \newtheorem{context}[theorem]{Context}

  \newtheorem*{defi}{Definition}
  \newtheorem*{acknowledgements}{Acknowledgements}



\begin{abstract}
We introduce a forcing technique to construct three-dimensional arrays of generic extensions through FS (finite support) iterations of ccc posets, which we refer to as 3D-coherent systems. We use them to produce models of new constellations in Cicho\'n's diagram, in particular, a model where the diagram can be separated into 7 different values. Furthermore, we show that this constellation of 7 values is consistent with the existence of a $\Delta^1_3$ well-order of the reals.
\end{abstract}
\maketitle

\section{Introduction}\label{SecIntro}

In this paper, we provide a generalization of the method of matrix iteration, to which we refer as \emph{3D-coherent systems} of iterations and which can be considered a natural extension of the matrix method to include a third dimension. That is, if a matrix iteration can be considered as a system of partial orders $\la \PP_{\alpha,\beta}:\alpha\leq\gamma, \beta\leq\delta\ra$ such that whenever $\alpha\leq\alpha'$ and $\beta\leq\beta'$ then $\PP_{\alpha,\beta}$ is a complete suborder of $\PP_{\alpha^\prime,\beta^\prime}$, then our 3D-coherent systems are systems of posets $\la \PP_{\alpha,\beta,\xi}:\alpha\leq\gamma,\beta\leq\delta,\xi\leq\pi\ra$ such that whenever $\alpha\leq\alpha^\prime$, $\beta\leq\beta^\prime$, $\xi\leq\xi^\prime$ then $\PP_{\alpha,\beta,\xi}$ is a complete suborder of $\PP_{\alpha^\prime,\beta^\prime,\xi^\prime}$. As an application of this method, we construct models where Cicho\'n's diagram is separated into different values, one of them with 7 different values. Moreover, these
models determine the value of $\afrak$, which is actually the same as the value of $\bfrak$, and we further show that such models can be produced so that they satisfy, additionally, the existence of a $\Delta^1_3$ well-order of the reals.

The method of matrix iterations, or {\emph{2D-coherent systems}} of iterations in our terminology, has already a long history. It was introduced by Blass and Shelah in~\cite{blassmatrix}, to show that consistently $\mathfrak{u}<\mathfrak{d}$, where $\mathfrak{u}$ is the ultrafilter number and $\mathfrak{d}$ is the dominating number. The method was further developed in~\cite{VFJB11}, where the terminology {\emph{matrix iteration}} appeared for the first time, to show that if $\kappa<\lambda$ are arbitrary regular uncountable cardinals then there is a generic extension in which $\mathfrak{a}=\mathfrak{b}=\kappa<\mathfrak{s}=\lambda$. Here $\mathfrak{a}$, $\mathfrak{b}$ and $\mathfrak{s}$ denote the almost disjointness, bounding and splitting numbers respectively. In~\cite{VFJB11}, the authors also introduce a new method for the preservation of a mad (maximal almost disjoint) family along
a matrix iteration, specifically a mad family added by $\Hor_\kappa$ (Hechler's poset for adding a mad family, see Definition \ref{DefHechlermad}), a method which is of particular importance for our current work. Later, classical preservation properties for matrix iterations were improved by Mej\'ia~\cite{mejia} to provide several examples of models where the cardinals in Cicho\'n's diagram assume many different values, in particular, a model with 6 different values. Since then, the question of {\emph{how many distinct values}} there can be simultaneously in Cicho\'n's diagram has been of interest for many authors, see for example~\cite{FGKS} (a model of 5 values concentrated on the right) and~\cite{GMS} (another model of 6 different values), and lies behind the development of many interesting forcing techniques. Very recently, the method of matrix iterations was used by Dow and Shelah~\cite{DowShelah} to solve a long-standing open
question in the area of cardinal characteristics of the continuum, namely, that it is consistent that the splitting number is singular.

Further motivation for this project was to determine the value of $\afrak$ in classical FS (finite support) iterations of ccc posets models where no dominating reals are added. To recall some examples, a classical result of Kunen~\cite{kunen80} states that, under CH, any Cohen poset preserves a mad family of the ground model. This result was improved by Steprans~\cite{steprans}, who showed that, after adding $\omega_1$-many Cohen reals, there is a mad family in the corresponding extension which is preserved by any further Cohen poset (without assuming CH). Additionally, Zhang~\cite{zhang} proved that, under CH, any finite support iteration of $\Eor$ (the standard poset adding an eventually different real, see Definition \ref{DefEDposet}) preserves a mad family from the ground model. As the family preserved in
Steprans' result is added by $\Cor_{\omega_1}=\Hor_{\omega_1}$, we considered the preservation theory of Brendle and the first author~\cite{VFJB11} to see in which cases a mad family added by $\Hor_\kappa$ (for an uncountable regular $\kappa$) can be preserved through FS iterations of ccc posets. If such an FS iteration can be redefined as a matrix iteration where $\Hor_\kappa$ is used to add a mad family as in \cite{VFJB11} and the preservation theory applies, then the mad family added by $\Hor_\kappa$ is preserved through the
iteration. Thanks to this and to the fact that random forcing and $\Eor$ fit into the preservation framework (Lemmas  \ref{randompresMAD} and~\ref{EpresMAD}), we generalize both Steprans' and Zhang's results by providing a general result about FS iterations preserving the mad family added by $\Hor_\kappa$ (Theorem
\ref{MainFSit}).

In view of the previous result, it is worth asking whether such a result can be extended to matrix iterations like those in~\cite{mejia}. By analogy, if it is possible to add an additional coordinate for $\Hor_\kappa$ to a matrix iteration and produce a 3D iteration (3D-coherent system in our notation) where the preservation
theory from~\cite{VFJB11} applies, then the mad family added by $\Hor_\kappa$ is preserved. Even more, the third dimension allows us to separate $\bfrak$ from other cardinals in Cicho\'n's diagram (which was not possible in~\cite{mejia}) and get a further division in Cicho\'n's diagram. In particular, the 3D-version of the matrix iteration from~\cite{mejia} for the consistency of 6 different values yields a model of 7 different values in Cicho\'n's diagram.

In addition, we show that these new constellations of Cicho\'n's diagram are consistent with the existence of a $\Delta^1_3$ well-order of the reals. Combinatorial properties of the real line (which can be expressed in terms of its cardinal characteristics) as well as the existence of nicely definable combinatorial objects (like maximal almost disjoint families) in the presence of a projective well-order on the reals have been investigated intensively in recent years.
In~\cite{VFSF} it is shown for example that various constellations involving $\mathfrak{a},\mathfrak{b}$ and $\mathfrak{s}$ are consistent with the existence of a $\Delta^1_3$ well-order, while in ~\cite{VFSFYK14} it is shown that every admissible assignment of $\aleph_1$ and $\aleph_2$ to the characteristics in Cicho\'n's diagram is consistent with the existence of such a projective well-order. There is one main distinction between the various known methods for generically adjoining projective well-orders: methods relying on countable support $S$-proper iterations like in~\cite{VFSF,VFSFYK14}, and methods using finite support iterations of ccc posets, e.g.~\cite{VFSFLZ11,VFSFAT12,VFSFLZ13}. In order to show that our new consistent constellations of Cicho\'n's diagram admit the existence of a $\Delta^1_3$ well-order of the reals, we further develop the second approach. Namely, we build up the method of almost disjoint coding which was introduced in~\cite{VFSFLZ11} and in particular answer one of the open
questions stated in~\cite{VFSFYK14}.

The paper is organized as follows. In Section \ref{SecPresProp} we present some well known preservation theorems. In Section \ref{SecCoherent} we introduce the notion of 3D-iteration and review the preservation properties for matrix iterations from~\cite{VFJB11,mejia} which can be applied quite directly to 3D-coherent systems (even to arbitrary coherent systems). In Section \ref{SecMad} we review the method of preservation of a mad family along a matrix iteration as introduced in~\cite{VFJB11} and obtain similar results regarding $\Eor$ and the random algebra. As a consequence, we prove in Theorem \ref{MainFSit} our generalization of Steprans' result discussed above, which is one of the main results of this paper.

Section \ref{SecMain} contains our main results about Cicho\'n's diagram. We evaluate the almost disjointness number in various constellations in which the value of $\mathfrak{a}$ was previously not known, and obtain a model in which there are 7 distinct values in Cicho\'n's diagram. Let $\theta_0\leq\theta_1\leq\kappa\leq\mu\leq\nu$ be regular uncountable cardinals, and let $\lambda\geq\nu$.

\begin{thm}
 Assume $\lambda^{<\theta_1}=\lambda$. Then, there is a ccc poset forcing $\add(\Nwf)=\theta_0$, $\cov(\Nwf)=\theta_1$, $\bfrak=\afrak=\kappa$, $\non(\Mwf)=\cov(\Mwf)=\mu$, $\dfrak=\nu$ and $\non(\Nwf)=\cfrak=\lambda$.
\end{thm}

Elaborating on the method of almost disjoint coding as developed in ~\cite{VFSFLZ11}, we show in Section \ref{SecDelta13} that the constellations of Section \ref{SecMain} are  consistent with the existence of a projective well-order of the reals whenever the associated cardinal values do not exceed $\aleph_\omega$ (even though we conjecture that the result remains true with arbitrarily large cardinal values). In particular, we outline the proof of the following:

\begin{thm}
In $L$, let $\theta_0<\theta_1<\kappa<\mu<\nu<\lambda$ be uncountable regular cardinals and, in addition, $\lambda <\aleph_\omega$. Then, there is a cardinal preserving forcing extension of $L$ in which there is a $\Delta^1_3$ well-order of the reals and, in addition, $\add(\Nwf)=\theta_0$, $\cov(\Nwf)=\theta_1$, $\bfrak=\afrak=\kappa$, $\non(\Mwf)=\cov(\Mwf)=\mu$, $\dfrak=\nu$ and $\non(\Nwf)=\cfrak=\lambda$.
\end{thm}

Section \ref{SecQ} contains some further discussions and open questions.

\bigskip

We recall some standard ccc posets we are going to use throughout this paper.

\begin{definition}[Standard forcing that adds an eventually different real]\label{DefEDposet}
   Define the forcing notion $\Eor$ with conditions of the form $(s,\varphi)$ where $s\in\omega^{<\omega}$ and $\varphi:\omega\to[\omega]^{<\aleph_0}$ such that $\exists{n<\omega}\forall{i<\omega}(|\varphi(i)|\leq n)$. Denote the minimal such $n$ by $\wdt(\varphi)$. The order in $\Eor$ is defined as $(t,\psi)\leq(s,\phi)$ iff $s\subseteq t$, $\forall{i<\omega}(\varphi(i)\subseteq\psi(i))$ and $\forall{i\in|t|\menos|s|}(t(i)\notin\varphi(i))$.
\end{definition}

Clearly $\Eor$ is $\sigma$-centered and adds a real which is eventually different from the reals in the ground model.  We will use also the following notation. If $\Omega$ is a non-empty set, $\Bor_\Omega$ is the cBa (complete Boolean algebra) $2^{\Omega\times\omega}/\Nwf(2^{\Omega\times\omega})$. Here, $\Nwf(2^{\Omega\times\omega})$ denotes the $\sigma$-ideal of null subsets of $2^{\Omega\times\omega}$ with respect to the standard product measure. Note that $\Bor_\Omega\simeq\Bor:=\Bor_\omega$ when $\Omega$ is countable. Also, for any non-empty set $\Gamma$, $\Bor_\Gamma:=\limdir\{\Bor_{\Omega}:\Omega\subseteq\Gamma\textrm{ countable}\}$. Denote by $\ral$ the class of all random algebras, that is, $\ral:=\{\Bor_\Gamma:\Gamma\neq\emptyset\}$. Recall Cohen forcing $\Cor_\Gamma:=\mathrm{Fn}(\Gamma\times\omega,2)$ which is the poset of finite partial functions from $\Gamma\times\omega$ to 2 ordered by reverse inclusion. Put $\Cor=\Cor_\omega$. Another well-known poset which we will make use of is the
localization poset (see for example~\cite{judabarto}). For convenience, we repeat its definition:

\begin{definition}\label{loc_poset} $\Loc$ is the poset of all $\varphi\in ([\omega]^{<\aleph_0})^\omega$ such that
\begin{enumerate}[(i)]
\item for all $n\in\omega$, $|\varphi(n)|\leq n$, and
\item there is a $k\in\omega$ such that for all but finitely many $n$, $|\varphi(n)|\leq k$.
\end{enumerate}
The extension relation is defined as follows: $\varphi' \leq \varphi$ if and only if $\varphi(n)\subseteq \varphi'(n)$ for all $n<\omega$.
\end{definition}

\section{Preservation properties for FS iterations}\label{SecPresProp}

We review the theory of preservation properties for FS iterations developed by Judah and Shelah \cite{jushe} and Brendle \cite{Br-Cichon}. A similar presentation also appears in \cite[Sect. 3]{GMS}.

\newcommand{\Rbf}{\mathbf{R}}

\begin{definition}\label{ContextFS}
   $\Rbf:=\la X,Y,\sqsubset\ra$ is a \emph{Polish relational system} if the following is satisfied:
   \begin{enumerate}[(i)]
      \item $X$ is a perfect Polish space,
      \item $Y$ is a non-empty analytic subspace of some Polish space and
      \item $\sqsubset=\bigcup_{n<\omega}\sqsubset_n$ for some increasing sequence $\la\sqsubset_n\ra_{n<\omega}$ of closed subsets of $X\times Y$ such that
            $(\sqsubset_n)^y=\{x\in X: x\sqsubset_n y\}$ is nwd (nowhere dense) for all $y\in Y$.
   \end{enumerate}
   For $x\in X$ and $y\in Y$, $x\sqsubset y$ is often read \emph{$y$ $\sqsubset$-dominates $x$}. A family $F\subseteq X$ is \emph{$\Rbf$-unbounded} if there is \underline{no} real in $Y$ that $\sqsubset$-dominates every member of $F$. Dually, $D\subseteq Y$ is a \emph{$\Rbf$-dominating} family if every member of $X$ is $\sqsubset$-dominated by some member of $D$. $\bfrak(\Rbf)$ denotes the least size of a $\Rbf$-unbounded family and $\dfrak(\Rbf)$ is the least size of a $\Rbf$-dominating family.

   Say that $x\in X$ is \emph{$\Rbf$-unbounded over a model $M$} if $x\not\sqsubset y$ for all $y\in Y\cap M$. Given a cardinal $\lambda$ say that $F\subseteq X$ is \emph{$\lambda$-$\Rbf$-unbounded} if, for any $Z\subseteq Y$ of size $<\lambda$, there is an $x\in F$ which is $\Rbf$-unbounded over $Z$.
\end{definition}

By (iii), $\la X,\Mwf(X),\in\ra$ is Tukey-Galois below $\Rbf$ where $\Mwf(X)$ denotes the $\sigma$-ideal of meager subsets of $X$. Therefore, $\bfrak(\Rbf)\leq\non(\Mwf)$ and $\cov(\Mwf)\leq\dfrak(\Rbf)$. Fix, for this section, a Polish relational system $\Rbf=\la X,Y,\sqsubset\ra$ and an uncountable regular cardinal $\theta$.

\begin{remark}\label{wlogPRS}
   Without loss of generality, $Y=\omega^\omega$ can be assumed. The reason is that, by the existence of a continuous surjection $f:\omega^\omega\to Y$, the Polish relational system $\Rbf':=\la X,\omega^\omega,\sqsubset'\ra$, where $x\sqsubset'_n z$ iff $x\sqsubset_n f(z)$, behaves much like $\Rbf$ in practice. Namely, $\Rbf$ is Tukey-Galois equivalent to $\Rbf'$ and moreover, the notions $\lambda$-$\Rbf$-unbounded and $\lambda$-$\Rbf'$-unbounded are equivalent. Also, for posets, the notions of $\theta$-$\Rbf$-good and $\theta$-$\Rbf'$-good (see the definition below) are equivalent.
\end{remark}

\begin{definition}[Judah and Shelah {\cite{jushe}}]\label{DefGood}
   A poset $\Por$ is \emph{$\theta$-$\Rbf$-good} if, for any $\Por$-name $\dot{h}$ for a real in $Y$, there is a non-empty $H\subseteq Y$ of size $<\theta$ such that $\Vdash x\not\sqsubset\dot{h}$ for any $x\in X$ that is $\Rbf$-unbounded over $H$.

   Say that $\Por$ is \emph{$\Rbf$-good} when it is $\aleph_1$-$\Rbf$-good.
\end{definition}

Definition \ref{DefGood} describes a property, respected by FS iterations, to preserve specific types of $\Rbf$-unbounded families. Concretely,
\begin{enumerate}[(a)]
  \item any $\theta$-$\Rbf$-good poset preserves every $\theta$-$\Rbf$-unbounded family from the ground model and
  \item FS iterations of $\theta$-cc $\theta$-$\Rbf$-good posets produce $\theta$-$\Rbf$-good posets.
\end{enumerate}
Posets that are $\theta$-$\Rbf$-good work to preserve $\bfrak(\Rbf)$ small and $\dfrak(\Rbf)$ large since, whenever $F$ is a $\theta$-$\Rbf$-unbounded family, $\bfrak(\Rbf)\leq|F|$ and $\theta\leq\dfrak(\Rbf)$.

Clearly, $\theta$-$\Rbf$-good implies $\theta'$-$\Rbf$-good whenever $\theta\leq\theta'$ and any poset completely embedded into a $\theta$-$\Rbf$-good poset is also $\theta$-$\Rbf$-good.

Consider the following particular cases of interest for our main results.

\begin{lemma}[{\cite[Lemma 4]{mejia}}]\label{smallGood}
   Any poset of size $<\theta$ is $\theta$-$\Rbf$-good. In particular, Cohen forcing is $\Rbf$-good.
\end{lemma}

\newcommand{\id}{\mathrm{id}}
\newcommand{\Ed}{\mathbf{Ed}}
\newcommand{\Dbf}{\mathbf{D}}
\newcommand{\Edb}{\mathbf{Ed}_b}
\newcommand{\Lc}{\mathbf{Lc}}

\begin{example}\label{ExmpGood}
   \begin{enumerate}[(1)]
      \item \emph{Preserving non-meager sets:} Consider the Polish relational system $\Ed:=\la\omega^\omega,\omega^\omega,\neq^*\ra$ where $x\neq^*y$ iff $x$ and $y$ are eventually different, that is, $x(i)\neq y(i)$ for all but finitely many $i<\omega$. By \cite[Thm. 2.4.1 and 2.4.7]{judabarto}, $\bfrak(\Ed)=\non(\Mwf)$ and $\dfrak(\Ed)=\cov(\Mwf)$.
      \item \emph{Preserving unbounded families:} Let $\Dbf:=\la\omega^\omega,\omega^\omega,\leq^*\ra$ be the Polish relational system where $x\leq^*y$ iff $x(i)\leq y(i)$ for all but finitely many $i<\omega$. Clearly, $\bfrak(\Dbf)=\bfrak$ and $\dfrak(\Dbf)=\dfrak$.

          Miller \cite{Mi} proved that $\Eor$ is $\Dbf$-good. Further,
          $\omega^\omega$-bounding posets,
like the random algebra, are $\Dbf$-good.

      \item \emph{Preserving null-covering families:} Let $b:\omega\to\omega\menos\{0\}$ such that $\sum_{i<\omega}\frac{1}{b(i)}<+\infty$ and let $\Edb:=\la\R_b,\R_b,\neq^*\ra$ be the Polish relational system where $\R_b:=\prod_{i<\omega}b(i)$. Since $\Edb$ is Tukey-Galois below $\la\Nwf(\R_b),\R_b,\not\ni\ra$ (for $x\in\R_b$ the set $\{y\in\R_b:\neg(x\neq^* y)\}$ has measure zero with respect to the standard Lebesgue measure on $\R_b$), $\cov(\Nwf)\leq\bfrak(\Edb)$ and $\dfrak(\Edb)\leq\non(\Nwf)$.

          By a similar argument as in \cite[Lemma $1^*$]{Br-Cichon}, any $\nu$-centered poset is $\theta$-$\Edb$-good for any $\nu<\theta$ infinite. In particular, $\sigma$-centered posets are $\Edb$-good.
      \item \emph{Preserving ``union of null sets is not null":} For each $k<\omega$ let $\id^k:\omega\to\omega$ such that $\id^k(i)=i^k$ for all $i<\omega$ and put $\Hwf:=\{\id^{k+1}:k<\omega\}$. Let $\Lc:=\la\omega^\omega,\Swf(\omega,\Hwf),\in^*\ra$ be the Polish relational system where
          \[\Swf(\omega,\Hwf):=\{\varphi:\omega\to[\omega]^{<\aleph_0}:\exists{h\in\Hwf}\forall{i<\omega}(|\varphi(i)|\leq h(i))\},\]
          and $x\in^*\varphi$ iff $\exists{n<\omega}\forall{i\geq n}(x(i)\in \varphi(i))$, which is read \emph{$x$ is localized by $\varphi$}. As a consequence of Bartoszy\'nski's characterization (see \cite[Thm. 2.3.9]{judabarto}), $\bfrak(\Lc)=\add(\Nwf)$ and $\dfrak(\Lc)=\cof(\Nwf)$.

          Any $\nu$-centered poset is $\theta$-$\Lc$-good for any $\nu<\theta$ infinite (see \cite{jushe}) so, in particular, $\sigma$-centered posets are $\Lc$-good. Moreover, subalgebras (not necessarily complete) of random forcing are $\Lc$-good as a consequence of a result of Kamburelis \cite{kamburelis}.
   \end{enumerate}
\end{example}

The following are the main general results concerning the preservation theory presented so far.

\begin{lemma}\label{CohenUnb}
   Let $\la\Por_\alpha\ra_{\alpha<\theta}$ be a $\lessdot$-increasing sequence of ccc posets and $\Por_\theta=\limdir_{\alpha<\theta}\Por_\alpha$. If $\Por_{\alpha+1}$ adds a Cohen real $\dot{c}_\alpha$ over $V^{\Por_\alpha}$ for any $\alpha<\theta$, then $\Por_\theta$ forces that $\{\dot{c}_\alpha:\alpha<\theta\}$ is a $\theta$-$\Rbf$-unbounded family of size $\theta$.
\end{lemma}

\begin{theorem}\label{MainPresThm}
   Let $\delta\geq\theta$ be an ordinal and $\la\Por_\alpha,\Qnm_\alpha\ra_{\alpha<\delta}$ an FS iteration of non-trivial $\theta$-$\Rbf$-good ccc posets. Then, $\Por_\delta$ forces $\bfrak(\Rbf)\leq\theta$ and $\dfrak(\Rbf)\geq|\delta|$.
\end{theorem}
\begin{proof}
  See \cite[Cor. 3.6]{GMS}.
\end{proof}

\section{Coherent systems of FS iterations}\label{SecCoherent}

\newcommand{\matit}{\mathbf{m}}
\newcommand{\tbf}{\mathbf{t}}

\begin{definition}[Relative embeddability]
   Let $M$ be a transitive model of ZFC (or a finite large fragment of it), $\Por\in M$ and $\Qor$ posets (the latter not necessarily in $M$). Say that \emph{$\Por$ is a complete subposet of $\Qor$ with respect to $M$}, denoted by $\Por\lessdot_M\Qor$, if $\Por$ is a suborder of $\Qor$ and every maximal antichain in $\Por$ that belongs to $M$ is also a maximal antichain in $\Qor$.
\end{definition}

Recall that in this case, if $N\supseteq M$ is another transitive model of ZFC with $\Qor\in N$ and $G$ is $\Qor$-generic over $N$ then $G\cap\Por$ is $\Por$-generic over $M$ and $M[G\cap\Por]\subseteq N[G]$. Moreover, if $\dot{\Por}'\in M$ is a $\Por$-name of a poset, $\dot{\Qor}'\in N$ is a $\Qor$-name of a poset and $\Vdash_{\Qor,N}\dot{\Por}'\lessdot_{M^{\Por}}\dot{\Qor}'$, then $\Por\ast\dot{\Por}'\lessdot_M\Qor\ast\dot{\Qor}'$. In particular, if $M=N=V$ (the universe), then $\Por\ast\dot{\Por}'\lessdot\Qor\ast\dot{\Qor}'$ whenever $\Por\lessdot\Qor$ and $\Vdash_{\Qor}\dot{\Por}'\lessdot_{V^\Por}\dot{\Qor}'$.

\newcommand{\sbf}{\mathbf{s}}

\begin{definition}[Coherent system of FS iterations]\label{DefCoherent}
   A \emph{coherent system (of FS iterations)} $\sbf$ is composed by the following objects:
   \begin{enumerate}[(I)]
     \item a partially ordered set $I^\sbf$ and an ordinal $\pi^\sbf$,
     \item a system of posets $\la\Por^\sbf_{i,\xi}:i\in I^\sbf,\xi\leq\pi^\sbf\ra$ such that
           \begin{enumerate}[(i)]
              \item $\Por^\sbf_{i,0}\lessdot\Por^\sbf_{j,0}$ whenever $i\leq j$ in $I^\sbf$, and
              \item $\Por^\sbf_{i,\eta}$ is the direct limit of $\la\Por^\sbf_{i,\xi}:\xi<\eta\ra$ for each limit $\eta\leq\pi^\sbf$,
           \end{enumerate}
     \item a sequence $\la\Qnm^\sbf_{i,\xi}:i\in I^\sbf,\xi<\pi^\sbf\ra$ where each $\Qnm^\sbf_{i,\xi}$ is a $\Por^\sbf_{i,\xi}$-name for a poset, $\Por^\sbf_{i,\xi+1}=\Por^\sbf_{i,\xi}\ast\Qnm^\sbf_{i,\xi}$ and $\Por^\sbf_{j,\xi}$ forces $\Qnm^\sbf_{i,\xi}\lessdot_{V^{\Por^\sbf_{i,\xi}}}\Qnm^\sbf_{j,\xi}$ whenever $i\leq j$ in $I^\sbf$ and $\Por^\sbf_{i,\xi}\lessdot\Por^\sbf_{j,\xi}$.
   \end{enumerate}
   Note that, for a fixed $i\in I^\sbf$, the posets $\la\Por^\sbf_{i,\xi}:\xi\leq\pi^\sbf\ra$ are generated by an FS iteration $\la\Por'_{i,\xi},\Qnm'_{i,\xi}:\xi<1+\pi^\sbf\ra$ where $\Qnm'_{i,0}=\Por^\sbf_{i,0}$ and $\Qnm'_{i,1+\xi}=\Qnm^\sbf_{i,\xi}$ for all $\xi<1+\pi^\sbf$. Therefore (by induction) $\Por'_{i,1+\xi}=\Por_{i,\xi}$ for all $\xi\leq\pi^\sbf$ and, thus, $\Por^\sbf_{i,\xi}\lessdot\Por^\sbf_{i,\eta}$ whenever $\xi\leq\eta\leq\pi^\sbf$.

   On the other hand, by Lemma  \ref{parallellimits}, $\Por^\sbf_{i,\xi}\lessdot\Por^\sbf_{j,\xi}$ whenever $i\leq j$ in $I^\sbf$ and $\xi\leq\pi^\sbf$.

   For $j\in I^\sbf$ and $\eta\leq\pi^\sbf$ we write $V^\sbf_{j,\eta}$ for the $\Por^\sbf_{j,\eta}$-generic extensions. Concretely, if $G$ is $\Por^\sbf_{j,\eta}$-generic over $V$, $V^\sbf_{j,\eta}:=V[G]$ and $V^\sbf_{i,\xi}:=V[\Por^\sbf_{i,\xi}\cap G]$ for all $i\leq j$ in $I^\sbf$ and $\xi\leq\eta$. Note that $V^\sbf_{i,\xi}\subseteq V^\sbf_{j,\eta}$.

   We say that the coherent system $\sbf$ has the \emph{ccc} if, additionally, $\Por^\sbf_{i,0}$ has the ccc and $\Por^\sbf_{i,\xi}$ forces that $\Qnm^\sbf_{i,\xi}$ has the ccc for each $i\in I^\sbf$ and $\xi<\pi^\sbf$. This implies that $\Por^\sbf_{i,\xi}$ has the ccc for all $i\in I^\sbf$ and $\xi\leq\pi^\sbf$.

   We consider the following particular cases.
   \begin{enumerate}[(1)]
      \item When $I^\sbf$ is a well-ordered set, we say that $\sbf$ is a \emph{2D-coherent system (of FS iterations)}.
      \item If $I^\sbf$ is of the form $\{i_0,i_1\}$ ordered as $i_0<i_1$, we say that $\sbf$ is a \emph{coherent pair (of FS iterations)}.
      \item If $I^\sbf=\gamma^\sbf\times\delta^\sbf$ where $\gamma^\sbf$ and $\delta^\sbf$ are ordinals and the order of $I^\sbf$ is defined as $(\alpha,\beta)\leq(\alpha',\beta')$ iff $\alpha\leq\alpha'$ and $\beta\leq\beta'$,
          we say that $\sbf$ is a \emph{3D-coherent system (of FS iterations).}
   \end{enumerate}

   For a coherent system $\sbf$ and a set $J\subseteq I^\sbf$, $\sbf|J$ denotes the coherent system with $I^{\sbf|J}=J$, $\pi^{\sbf|J}=\pi^\sbf$ and the posets and names corresponding to (II) and (III) defined as for $\sbf$. And if $\eta\leq\pi^\sbf$, $\sbf\frestr\eta$ denotes the coherent system with $I^{\sbf\upharpoonright\eta}=I^\sbf$, $\pi^{\sbf\upharpoonright\eta}=\eta$ and the posets for (II) and (III) defined as for $\sbf$. Note that, if $i_0<i_1$ in $I^\sbf$, then $\sbf|\{i_0,i_1\}$ is a coherent pair and $\sbf|\{i_0\}$ corresponds just to the FS iteration $\la\Por'_{i_0,\xi},\Qnm'_{i_0,\xi}:\xi<1+\pi^\sbf\ra$ (see the comment after (III)).

   If $\tbf$ is a 3D-coherent system, for $\alpha<\gamma^\tbf$, $\tbf_\alpha:=\tbf|\{(\alpha,\beta):\beta<\delta^\tbf\}$ is a 2D-coherent system where $I^{\tbf_\alpha}$ has order type $\delta^\tbf$. For $\beta<\delta^\tbf$, $\tbf^\beta:=\tbf|\{(\alpha,\beta):\alpha<\delta^\tbf\}$ is a 2D-coherent system where $I^{\tbf^\beta}$ has order type $\gamma^\tbf$.

   In particular, the upper indices $\sbf$ are omitted when there is
   no risk of ambiguity.

\end{definition}

Concerning consistency results about cardinal characteristics of the real line, Blass and Shelah \cite{blassmatrix} produced the first 2D-coherent system to obtain that $\ufrak<\dfrak$ is consistent with large continuum. This was followed by new consistency results by Brendle and Fischer \cite{VFJB11} and Mej\'ia \cite{mejia} where Blass' and Shelah's construction (which consists, basically, of 2D-coherent systems as formalized in Definition \ref{DefCoherent}(1)) is formulated and improved. For their results, the main features of the produced matrix of generic extensions $\la V_{\alpha,\xi}:\alpha\leq\gamma,\xi\leq\pi\ra$ from a 2D-coherent system $\mathbf{m}$, as illustrated in Figure \ref{FigMatit}, are:

\begin{figure}
\begin{center}
     \includegraphics[scale=0.72]{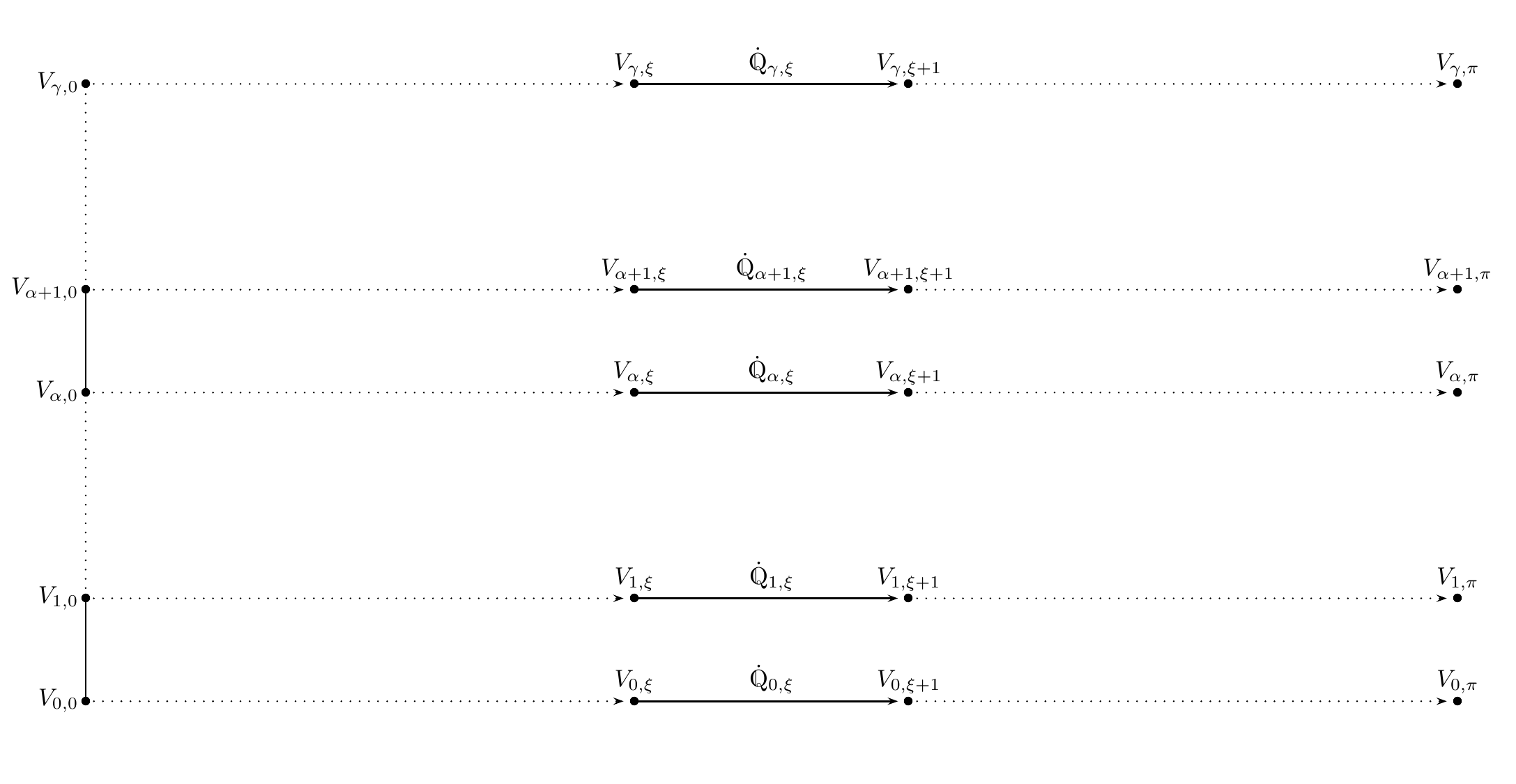}
     \caption{Matrix of generic extensions (2D-coherent system).}
     \label{FigMatit}
\end{center}
\end{figure}

\begin{enumerate}[(F1)]
   \item For $\alpha<\gamma$, there is a real $c_\alpha\in V_{\alpha+1,0}$ which ``diagonalizes" $V_{\alpha,0}$ (e.g., $\Rbf$-unbounded over $V_{\alpha,0}$ for a fixed Polish relational system $\Rbf$, or diagonalizes it in the sense of Definition \ref{DefpresDiag}) and, through the coherent pair $\mathbf{m}|\{\alpha,\alpha+1\}$, $c_\alpha$ also diagonalizes all the models in the $\alpha$-th row, that is, $V_{\alpha,\xi}$ for all $\xi\leq\pi$ (Lemmas \ref{parallellimits} and \ref{limitMAD}).
   \item Assume that $\gamma$ (the top level of the matrix) has uncountable cofinality. Given any column of the matrix, any real in the model of the top is actually in some of the models below, that is, $\R\cap V_{\gamma,\xi}=\bigcup_{\alpha<\gamma}\R\cap V_{\alpha,\xi}$ for every $\xi\leq\pi$ (Lemma \ref{GenNoNewReals} and Corollary \ref{Nonewreals}).
\end{enumerate}

\begin{figure}
\begin{center}
     \includegraphics{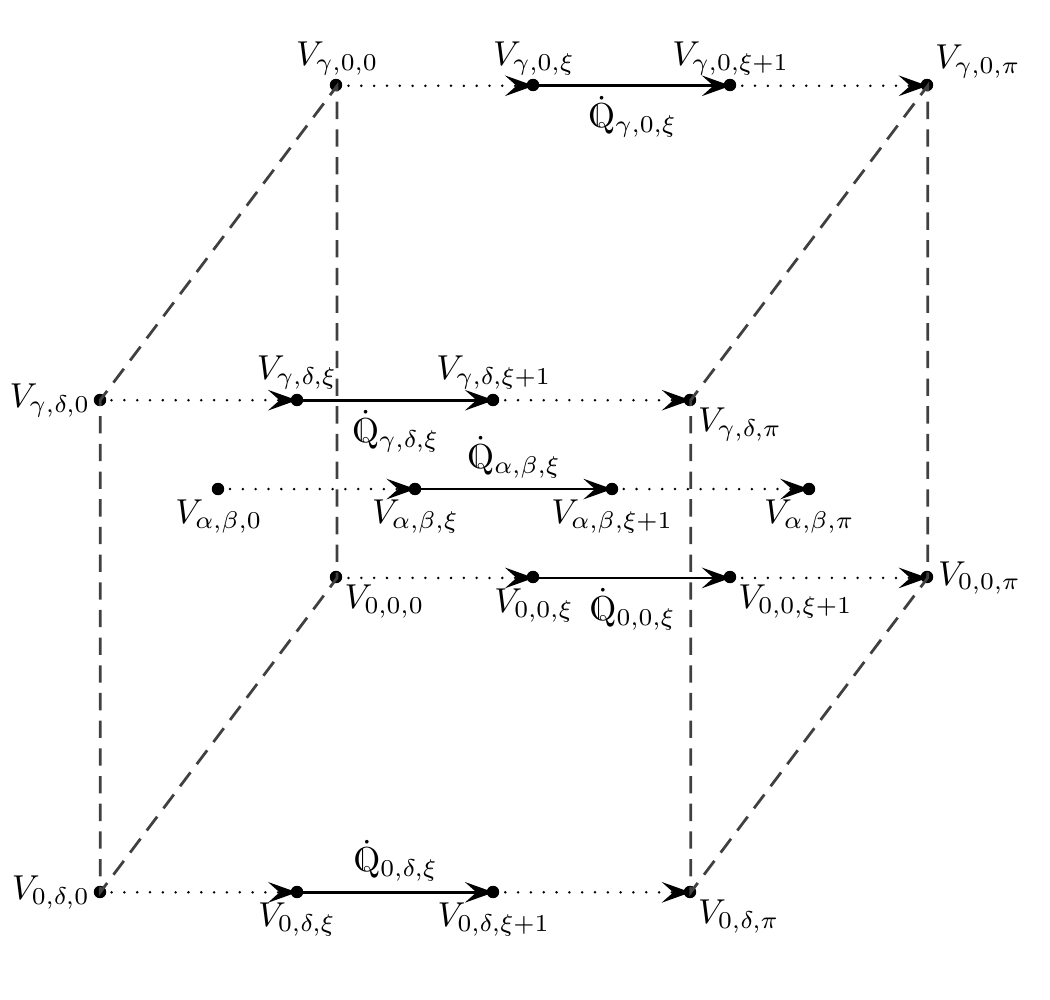}
     \caption{3D rectangle of generic extensions (3D-coherent system).}
     \label{Figcubeit}
\end{center}
\end{figure}

To prove the main results of this paper, we extend this approach to 3D rectangles of generic extensions which help us separate more cardinal invariants at the same time. In a similar fashion as a matrix above, such a construction starts with a matrix of posets and ``coherent" FS iterations emanate from each poset, which is formalized in Definition \ref{DefCoherent}(3) as 3D-coherent systems. Figure \ref{Figcubeit} illustrates this idea. More generally, Definition \ref{DefCoherent} can be used to define multidimensional rectangles of generic extensions, though applications are unknown for dimensions $\geq 4$.

The feature (F1) can also be applied in general to coherent systems of
FS iterations since any such system is composed of several coherent
pairs of FS iterations. For coherent pairs, (F1) for
``$\Rbf$-unbounded over a model" has been well understood in \cite{blassmatrix,VFJB11,mejia} whose results we review below. For the remainder of this section, fix $M\subseteq N$ transitive models of ZFC and a Polish relational system $\Rbf=\la X,Y,\sqsubset\ra$ coded in $M$ (in the sense that all its components are coded in $M$).

Recall that $\Sor$ is a \emph{Suslin ccc poset} if it is a
$\boldsymbol{\Sigma}^1_1$ subset of $\omega^\omega$ (or another
uncountable Polish space) and both its order and
incompatibility relations are $\boldsymbol{\Sigma}^1_1$. Note that if
$\Sor$ is coded in $M$ then $\Sor^M\lessdot_M\Sor^N$.

\begin{lemma}[{\cite[Thm. 7]{mejia}}]\label{Suslingoodparallel}
   Let $\Sor$ be a Suslin ccc poset coded in $M$. If $M\models``\Sor$ is $\Rbf$-good" then, in $N$, $\Sor^N$ forces that every real in $X\cap N$ which is $\Rbf$-unbounded over $M$ is $\Rbf$-unbounded over $M^{\Sor^M}$.
\end{lemma}

\begin{corollary}
   Let $\Gamma\in M$ be a non-empty set. If $M\models``\Bor_\Gamma$ is $\Rbf$-good" then $\Bor_\Gamma^N$, in $N$, forces that every real in $X\cap N$ which is $\Rbf$-unbounded over $M$ is $\Rbf$-unbounded over $M^{\Bor_\Gamma^M}$.
\end{corollary}

\begin{lemma}[{\cite[Lemma 11]{VFJB11}}, see also {\cite[Lemma 5.13]{mejia-temp}}]\label{Fixedparallel}
   Assume $\Por\in M$ is a poset. Then, in $N$, $\Por$ forces that every real in $X\cap N$ which is $\Rbf$-unbounded over $M$ is $\Rbf$-unbounded over $M^{\Por}$.
\end{lemma}

\begin{lemma}[Blass and Shelah {\cite{blassmatrix}}, {\cite[Lemmas 10, 12 and 13]{VFJB11}}]\label{parallellimits}
   Let $\sbf$ be a coherent pair of FS iterations as in Definition \ref{DefCoherent}(2). Then, $\Por_{i_0,\xi}\lessdot\Por_{i_1,\xi}$ for all $\xi\leq\pi$.

   Moreover, if $\dot{c}$ is a $\Por_{i_1,0}$-name of a real in $X$, $\pi$ is limit and $\Por_{i_1,\xi}$ forces that $\dot{c}$ is $\Rbf$-unbounded over $V_{i_0,\xi}$ for all $\xi<\pi$, then $\Por_{i_1,\pi}$ forces that $\dot{c}$ is $\Rbf$-unbounded over $V_{i_0,\pi}$.
\end{lemma}

Note that if $c$ is a Cohen real over $M$ then $c$ is $\Rbf$-unbounded over $M$ by Definition \ref{ContextFS}(iii). In fact, all the unbounded reals used in our applications are actually Cohen.

Now we turn to discuss feature (F2). We aim to have such a property
for 3D-coherent systems but, as they are composed of several
2D-coherent systems, it is enough to understand (F2) for 2D-coherent
systems. This was already noted in \cite{blassmatrix} and
formalized in \cite[Lemma 15]{VFJB11} (see Corollary
\ref{Nonewreals}), which we generalize as follows.

\begin{lemma}\label{GenNoNewReals}
   Let $\matit$ be a ccc 2D-coherent system with $I^\matit=\gamma+1$ an ordinal and $\pi^\matit=\pi$. Assume that
   \begin{enumerate}[(i)]
     \item $\gamma$ has uncountable cofinality,
     \item $\Por_{\gamma,0}$ is the direct limit of $\la\Por_{\alpha,0}:\alpha<\gamma\ra$, and
     \item for any $\xi<\pi$, $\Por_{\gamma,\xi}$ forces ``$\Qnm_{\gamma,\xi}=\bigcup_{\alpha<\gamma}\Qnm_{\alpha,\xi}$" whenever $\Por_{\gamma,\xi}$ is the direct limit of $\la\Por_{\alpha,\xi}:\alpha<\gamma\ra$.
   \end{enumerate}
   Then, for any $\xi\leq\pi$, $\Por_{\gamma,\xi}$ is the direct limit of $\la\Por_{\alpha,\xi}:\alpha<\gamma\ra$. In particular, $\Por_{\gamma,\xi}$ forces that $\R\cap V_{\gamma,\xi}=\bigcup_{\alpha<\gamma}\R\cap V_{\alpha,\xi}$.
\end{lemma}
\begin{proof}
   We proceed by induction on $\xi$. The case when $\xi$ is not successor is clear, so we just need to deal with the successor step. Assume that the conclusion holds for $\xi$. If $p\in\Por_{\gamma,\xi+1}$ then $p=(r,\dot{q})$ where $r\in\Por_{\gamma,\xi}$ and $\dot{q}$ is a $\Por_{\gamma,\xi}$-name of a member of $\Qnm_{\gamma,\xi}$. By (iii), there is a maximal antichain $\{p_n:n<\omega\}$ in $\Por_{\gamma,\xi}$ such that $p_n$ decides $\dot{q}=\dot{q}_n\in\Qnm_{\alpha_n,\xi}$ for some $\alpha_n<\gamma$ and some $\Por_{\alpha_n,\xi}$-name $\dot{q}_n$.\footnote{It is implicit in this proof that the names considered for the members of $\Qnm_{\alpha,\gamma}$ are \emph{canonical} in the sense described by the mentioned maximal antichains.} By (i), (ii) and the induction hypothesis, there is an $\alpha<\gamma$ above all $\alpha_n$ such that $\{p_n:n<\omega\}\subseteq\Por_{\alpha,\xi}$ and $r\in\Por_{\alpha,\xi}$. Therefore, $\dot{q}$ is a $\Por_{\alpha,\xi}$-name of a member of $\Qnm_{\alpha,\xi}$ and $p\in\Por_
{\alpha,\xi+1}$.
\end{proof}

The 2D and 3D-coherent systems constructed to prove our main results can be classified in terms of the following notion.

\begin{definition}[Standard coherent system of FS iterations]\label{DefStandMatit}
   A ccc coherent system of FS iterations $\sbf$ is \emph{standard} if
   \begin{enumerate}[(I)]
     \item it consists, additionally, of:
   \begin{enumerate}[(i)]
      \item a partition $\la S^{\sbf},C^{\sbf}\ra$ of $\pi^{\sbf}$,
      \item a function $\Delta^{\sbf}:C^{\sbf}\to I^{\sbf}$ so that $\Delta^{\sbf}(i)$ is not maximal in $I^\sbf$ for all $i\in C^\sbf$,
      \item a sequence $\la\Sor^{\sbf}_\xi:\xi\in S^{\sbf}\ra$ where each $\Sor^{\sbf}_\xi$ is either a Suslin ccc poset or a random algebra, and
      \item a sequence $\la\Qnm^{\sbf}_\xi:\xi\in C^{\sbf}\ra$ such that each $\Qnm^{\sbf}_\xi$ is a $\Por^{\sbf}_{\Delta^\sbf(\xi),\xi}$-name of a poset which is forced to be ccc by $\Por^{\sbf}_{i,\xi}$ for all $i\geq\Delta^\sbf(\xi)$ in $I^\sbf$, and
   \end{enumerate}
   \item it satisfies, for any $i\in I^\sbf$ and $\xi<\pi^{\sbf}$, that
         \[\Qnm^{\sbf}_{i,\xi}=\left\{\begin{array}{ll}
              (\Sor^{\sbf}_\xi)^{V^{\sbf}_{i,\xi}} & \textrm{if $\xi\in S^{\sbf}$}\\
              \Qnm^{\sbf}_\xi & \textrm{if $\xi\in C^{\sbf}$ and $i\geq\Delta^{\sbf}(\xi)$,}\\
              \mathds{1} & \textrm{otherwise.}
          \end{array}\right.\]
   \end{enumerate}
   As in Definition \ref{DefCoherent}, the upper index $\sbf$ may be omitted when it is clear from the context.
\end{definition}

All the standard coherent systems in this paper are constructed by recursion on $\xi<\pi$. To be more precise, we start with some partial order of ccc posets $\la \Por_{i,0}:i\in I\ra$ as in Definition \ref{DefCoherent}(II)(i), fix the partition in (I)(i) and, by recursion, the posets $\Por_{i,\xi}$ and names $\Qnm_{i,\xi}$ for all $i\in I$, along with the function $\Delta$ and the sequence of Suslin ccc posets in (I)(iii) (though in some cases $\Delta$ and the sequence of Suslin ccc posets are fixed before the recursion), are defined as follows: when $\Por_{i,\xi}$ has been constructed for all $i\in I$, we distinguish the cases $\xi\in S$ and $\xi\in C$. In the first case, $\Sor_\xi$ is chosen; in the second, we choose $\Delta(\xi)$ and then we define the ($\Por_{\Delta(\xi),\xi}$-name of a) poset $\Qnm_\xi$ as in (I)(iv). After this, the iterations continue with $\Por_{i,\xi+1}=\Por_{i,\xi}\ast\Qnm_{i,\xi}$ as indicated in (II). It is clear that the requirements in Definition \ref{DefCoherent} for a ccc
coherent system are satisfied.

In practice, a standard coherent system as above is constructed by using posets adding generic reals and the cases whether $\xi\in S$ or $\xi\in C$ indicate how generic the real is. Namely, when $\xi\in S$, $\Sor_\xi$ adds a real that is generic over $V_{i,\xi}$ for all $i\in I$, which means that we add a \emph{full generic real} at stage $\xi$; on the other hand, when $\xi\in C$ we just add a \emph{restricted generic} in the sense that $\Qnm_\xi$ adds a real which is generic over $V_{\Delta(\xi),\xi}$ but \underline{not} necessarily over $V_{i,\xi}$ when $i\not\leq\Delta(\xi)$, for instance, if $\Qnm_\xi$ is a name for $\Dor^{V_{\Delta(\xi),\xi}}$, at the $\xi$-step the Hechler real added is generic only over $V_{\Delta(\xi),\xi}$. This approach of adding full and restricted generic reals is useful for controlling many cardinal invariants at the same time like in \cite{blassmatrix,VFJB11,mejia} and this work.

It is clear that any standard 2D-coherent system satisfies the hypothesis (iii) of Lemma \ref{GenNoNewReals} whenever (i) and (ii) are satisfied. Therefore,

\begin{corollary}[{\cite[Lemma 15]{VFJB11}}]\label{Nonewreals}
   If $\matit$ is a standard 2D-coherent system with $I^\matit=\gamma+1$ and an ordinal and $\pi^\matit=\pi$ satisfying (i) and (ii) of Lemma \ref{GenNoNewReals} then, for any $\xi\leq\pi$, $\Por_{\gamma,\xi}$ is the direct limit of $\la\Por_{\alpha,\xi}:\alpha<\gamma\ra$. In particular, $\Por_{\gamma,\xi}$ forces that $\R\cap V_{\gamma,\xi}=\bigcup_{\alpha<\gamma}\R\cap V_{\alpha,\xi}$.
\end{corollary}

The results presented in this section can be summarized in the following result.

\begin{theorem}[{\cite[Thm. 10 \& Cor. 1]{mejia}}]\label{MainMatit}
   Let $\matit$ be a standard 2D-coherent system with $I^\matit=\gamma+1$ (an ordinal), $\pi^\matit=\pi$ and $\Rbf=\la X,Y,\sqsubset\ra$ a Polish relational system coded in $V$. Assume that
   \begin{enumerate}[(i)]
      \item for any $\xi\in S$ and $\alpha\leq\gamma$, $\Por_{\alpha,\xi}$ forces that $\Qnm_{\alpha,\xi}=\Sor_{\xi}^{V_{\alpha,\xi}}$ is $\Rbf$-good and
      \item for any $\alpha<\gamma$ there is a $\Por_{\alpha+1,0}$-name $\dot{c}_\alpha$ of a $\Rbf$-unbounded member of $X$ over $V_{\alpha,0}$.
   \end{enumerate}
   Then, for any $\xi\leq\pi$ and $\alpha<\gamma$, $\Por_{\alpha+1,\xi}$ forces that $\dot{c}_\alpha$ is $\Rbf$-unbounded over $V_{\alpha,\xi}$. In addition, if $\matit$ satisfies (i) and (ii) of Lemma \ref{GenNoNewReals} then $\Por_{\gamma,\pi}$ forces $\bfrak(\Rbf)\leq\cf(\gamma)\leq\dfrak(\Rbf)$.
\end{theorem}
\begin{proof}
   The first statement is a direct consequence of Lemmas \ref{Suslingoodparallel}, \ref{Fixedparallel} and \ref{parallellimits}. For the second statement, note that Corollary \ref{Nonewreals} implies that, in $V_{\gamma,\pi}$, $\{c_{\alpha_\eta}:\eta<\cf(\gamma)\}$ is a $\cf(\gamma)$-$\Rbf$-unbounded family where $\la\alpha_\eta:\eta<\cf(\gamma)\ra\in V$ is an increasing cofinal sequence of $\gamma$, so $\bfrak(\Rbf)\leq\cf(\gamma)\leq\dfrak(\Rbf)$ follows.
\end{proof}

\section{Preservation of Hechler mad families}\label{SecMad}

We review from \cite{VFJB11} the theory of preserving, through coherent pairs of FS iterations, a mad family added by Hechler's poset for adding an a.d. family (see Definition \ref{DefHechlermad}). This theory is quite similar to the approach in Section \ref{SecCoherent}. Additionally, we show in Lemmas \ref{EpresMAD} and \ref{randompresMAD} that random forcing $\Bor$ and the eventually different forcing $\Eor$ fit well into this framework.

\begin{definition}[Hechler {\cite{Hechlermad}}]\label{DefHechlermad}
  For a set $\Omega$ define the poset $\Hor_\Omega:=\{p:F_p\times n_p\to 2:F_p\in[\Omega]^{<\aleph_0}\textrm{\ and }n_p<\omega\}$. The order is given by $q\leq p$ iff $p\subseteq q$ and, for any $i\in n_q\menos n_p$, there is \underline{at most one} $z\in F_p$ such that $q(z,i)=1$.
\end{definition}

If $G$ is $\Hor_\Omega$-generic over $V$ then $A=A_G:=\{a_z:z\in\Omega\}$ is an a.d. family where $a_z\subseteq\omega$ is defined as $i\in a_z$ iff $p(z,i)=1$ for some $p\in G$. Moreover, $V[G]=V[A]$ and, when $\Omega$ is uncountable, $A$ is mad in $V[G]$ (see \cite{Hechlermad}).

If $\Omega\subseteq\Omega'$ it is clear that
$\Hor_\Omega\lessdot\Hor_{\Omega'}$ and even the ($\Hor_\Omega$-name
of the) quotient $\Hor_{\Omega'}/\Hor_{\Omega}$ is nicely expressed
(see, e.g., \cite[\S2]{VFJB11}). On the other hand, if $\Cwf$ is an
$\subseteq$-chain of sets then
$\Hor_{\bigcup\Cwf}=\limdir_{\Omega\in\Cwf}\Hor_\Omega$. Therefore, if
$\gamma$ is an ordinal, $\Hor_\gamma$ can be obtained by an FS
iteration of length $\gamma$ where $\Hor_\alpha$ is the poset obtained
in the $\alpha$-th stage of the iteration and
$\Hor_{\alpha+1}/\Hor_\alpha$, which is $\sigma$-centered, is the
$\alpha$-th iterand. Since $\Hor_\Omega$ only depends on the size of
$\Omega$, this implies that $\Hor_\Omega$ has precaliber $\omega_1$
(though this can be proved directly by a $\Delta$-system
argument). Moreover, if $\Omega$ is non-empty and countable then
$\Hor_\Omega\simeq\Cor$ and, if $|\Omega|=\aleph_1$, then
$\Hor_\Omega\simeq\Cor_{\omega_1}$.

From now on, fix transitive models $M\subseteq N$ of ZFC. We define below a diagonalization property to preserve mad families like the one added by Hechler's poset.

\begin{definition}[{\cite[Def. 2]{VFJB11}}]\label{DefpresDiag}
   Let $A=\la a_z\ra_{z\in\Omega}\in M$ be a family of infinite subsets of $\omega$ and $a^*\in[\omega]^{\aleph_0}$ (not necessarily in $M$). Say that \emph{$a^*$ diagonalizes $M$ outside $A$} if, for all $h\in M$, $h:\omega\times[\Omega]^{<\aleph_0}\to\omega$ and for any $m<\omega$, there are $i\geq m$ and $F\in[\Omega]^{<\aleph_0}$ such that $[i,h(i,F))\menos\bigcup_{z\in F}a_z\subseteq a^*$.
\end{definition}

Given a collection $A$ of subsets of $\omega$, \emph{the ideal generated by $A$} is defined as

\[\Iwf(A):=\{x\subseteq\omega:x\subseteq^*\bigcup_{a\in F}a\textrm{\ for some finite }F\subseteq A\}.\]

\begin{lemma}[{\cite[Lemma 3]{VFJB11}}]\label{DiagMain}
   If $a^*$ diagonalizes $M$ outside $A$ then $|a^*\cap x|=\aleph_0$ for any $x\in M\menos\Iwf(A)$.
\end{lemma}

\begin{corollary}\label{DiagMad}
   Let $\gamma$ be an ordinal of uncountable cofinality and let $\langle M_\alpha\rangle_{\alpha\leq\gamma}$ be an increasing sequence of transitive ZFC models such that $[\omega]^{\aleph_0}\cap M_\gamma=\bigcup_{\alpha<\gamma}[\omega]^{\aleph_0}\cap M_\alpha$. Assume that $A=\{a_\alpha:\alpha<\gamma\}\in M_\gamma$ is a family of infinite subsets of $\omega$ such that, for any $\alpha<\gamma$, $A\frestr\alpha\in M_\alpha$ and $a_\alpha\in M_{\alpha+1}$ diagonalizes $M_\alpha$ outside $A\frestr\alpha$. \underline{Then}, for any $x\in[\omega]^{\aleph_0}\cap M_\gamma$, there exists an $\alpha<\gamma$ such that $|x\cap a_\alpha|=\aleph_0$. If, additionally, $A$ is almost disjoint, then $A$ is mad in $M_\gamma$.
\end{corollary}

\begin{lemma}[{\cite[Lemma 4]{VFJB11}}]\label{HechlerMad}
   Let $\Omega$ be a set, $z^*\in\Omega$ and $A:=\{a_z:z\in\Omega\}$ the a.d. family added by $\Hor_\Omega$. Then, $\Hor_\Omega$ forces that $a_{z^*}$ diagonalizes $V^{\Hor_{\Omega\menos\{z^*\}}}$ outside $A\frestr(\Omega\menos\{z^*\})$
\end{lemma}

Though it is well-known that, for $\Omega$ uncountable, the a.d. family added by $\Hor_\Omega$ is mad (as mentioned earlier), this follows from Corollary \ref{DiagMad} and Lemma \ref{HechlerMad} since $\Hor_\Omega\cong\Hor_\gamma$ for some ordinal $\gamma$ of uncountable cofinality.

The main idea for mad preservation in \cite{VFJB11} is that, when ccc 2D-coherent systems are constructed, the first column, along with a mad family $A=\{a_\alpha:\alpha<\gamma\}$, satisfies the hypothesis of Corollary \ref{DiagMad} (e.g. $\Por_{\alpha,0}=\Hor_\alpha$ for all $\alpha\leq\gamma$) and each $a_\alpha$ is preserved to diagonalize the models in the $\alpha$-th row outside $A\frestr\alpha$ (that is, the second case of (F1) at the beginning of Section \ref{SecCoherent}). For this purpose, we present the following results related to the preservation of the property in Definition \ref{DefpresDiag} through coherent pairs of iterations.

\begin{lemma}[{\cite[Lemma 11]{VFJB11}}]\label{FixedDiag}
   Let $\Por\in M$ be a poset. If $N\models``a^*$ diagonalizes $M$ outside $A$" then
   \[N^\Por\models``\textrm{$a^*$ diagonalizes $M^\Por$ outside $A$}".\]
\end{lemma}

\begin{corollary}\label{CohenDiag}
   If $N\models``a^*$ diagonalizes $M$ outside $A$" then
   \[N^{\Cor^N}\models``\textrm{$a^*$ diagonalizes $M^{\Cor^M}$ outside $A$}".\]
\end{corollary}

\begin{lemma}\label{EpresMAD}
   If $N\models``a^*$ diagonalizes $M$ outside $A$" then
   \[N^{\Eor^N}\models``\textrm{$a^*$ diagonalizes $M^{\Eor^M}$ outside $A$}".\]
\end{lemma}
\begin{proof}
   Let $\dot{h}\in M$ be an $\Eor$-name for a function from $\omega\times[\Omega]^{<\aleph_0}$ into $\omega$. Work within $M$ and fix a non-principal ultrafilter $D$ on $\omega$ (in $M$). For $s\in\omega^{<\omega}$ and $n<\omega$ define $h_{s,n}:\omega\times[\Omega]^{<\aleph_0}\to\omega+1$ as
   \[h_{s,n}(i,F)=\min\{j<\omega:(\forall{\varphi,\ \wdt(\varphi)\leq n})((s,\varphi)\nVdash\dot{h}(i,F)>j)\}.\]
   \begin{clm}
      $h_{s,n}(i,F)\in\omega$ for all $i<\omega$ and $F\in[\Omega]^{<\aleph_0}$.
   \end{clm}
   \begin{proof}
      Assume not, so there is a sequence of slaloms $\la\varphi_j\ra_{j<\omega}$ of width $\leq n$ such that $(s,\varphi_j)\Vdash\dot{h}(i,F)>j$. Define the slalom $\varphi^*$ as
      \[\varphi^*(i)=\{m<\omega:\{j<\omega:m\in\varphi_j(i)\}\in D\}.\]
      Since $D$ is a filter, $\wdt(\varphi^*)\leq n$, so $(s,\varphi^*)\in\Eor$. Now, there are $(t,\psi)\leq(s,\varphi^*)$ and $j_0<\omega$ such that $(t,\psi)\Vdash\dot{h}(i,F)=j_0$. By the definition of $\varphi^*$ and since $D$ is an ultrafilter,
      \[\{j<\omega:\forall{i\in|t|\menos|s|}(t(i)\notin\varphi_j(i))\}\in D\]
      so that set is infinite. For any $j>j_0$ in that set, $(t,\psi)$ is compatible with $(s,\varphi_j)$ and, therefore, any common stronger condition forces $j_0=\dot{h}(i,F)>j$, a contradiction.
   \end{proof}
   Now, in $N$, fix $m<\omega$ and $p=(s,\varphi)\in\Eor^N$ with $n:=\wdt(\varphi)$. As $a^*$ diagonalizes $M$ outside $A$, there are $i\geq m$ and $F\in[\Omega]^{<\aleph_0}$ such that $[i,h_{s,n}(i,F))\menos\bigcup_{z\in F}a_z\subseteq a^*$. By definition of $h_{s,n}$, $(\forall{\varphi,\ \wdt(\varphi)\leq n})((s,\varphi)\nVdash\dot{h}(i,F)>h_{s,n}(i,F))$ is a true $\boldsymbol{\Pi}^1_1$-statement in $M$ so, by absoluteness, it is also true in $N$. Therefore, there is a $q\in\Eor^N$ stronger than $p$ that forces $\dot{h}(i,F)\leq h_{s,n}(i,F)$ and then we conclude that $q$ forces $[i,\dot{h}(i,F))\menos\bigcup_{z\in F}a_z\subseteq a^*$.
\end{proof}

\begin{lemma}\label{randompresMAD}
   If $N\models``a^*$ diagonalizes $M$ outside $A$" then
   \[N^{\Bor^N}\models``\textrm{$a^*$ diagonalizes $M^{\Bor^M}$ outside $A$}".\]
\end{lemma}
\begin{proof}
   In the standard proof that $\Bor$ is $\omega^\omega$-bounding (see for example~\cite{judabarto}) it is shown that, for any $p\in\Bor$, $\epsilon\in(0,1)$ and $\dot{x}$ a $\Bor$-name for a real in $\omega^\omega$, there are $q\leq p$ and $g\in\omega^\omega$ such that $q\Vdash\dot{x}\leq g$ and $\lambda(p\menos q)\leq\epsilon\lambda(p)$ where $\lambda$ is the Lebesgue measure. We are going to use this fact to prove the lemma.

   Fix $\dot{h}\in M$ a $\Bor$-name for a function from $\omega\times[\Omega]^{<\aleph_0}$ to $\omega$, $p\in\Bor^N$ and $m<\omega$. By the Lebesgue density Theorem there is a clopen non-empty set $C$ such that $\lambda(C\menos p)<\frac{1}{4}\lambda(C)$. Now, in $M$, find $g:\omega\times[\Omega]^{<\aleph_0}\to\omega$ such that, for any $F\in[\Omega]^{<\aleph_0}$, there is a $q_F\leq C$ in $\Bor$ with $\lambda(C\menos q_F)\leq\frac{1}{4}\lambda(C)$ that forces $\forall{i<\omega}(\dot{h}(i,F)\leq g(i,F))$. Then, in $N$, there are $i\geq m$ and $F\in[\Omega]^{<\aleph_0}$ such that $[i,g(i,F))\menos\bigcup_{z\in F}a_z\subseteq a^*$, so $q_F$ forces $[i,\dot{h}(i,F))\menos\bigcup_{z\in F}a_z\subseteq a^*$. As $\lambda(p\cap q_F)>\frac{1}{2}\mu(C)$, $p\cap q_F\in\Bor^N$ is stronger than $p$ and forces $[i,\dot{h}(i,F))\menos\bigcup_{z\in F}a_z\subseteq a^*$.
\end{proof}

\begin{corollary}\label{randomalgpresMAD}
   Let $\Gamma\in M$ be a non-empty set. If $N\models``a^*$ diagonalizes $M$ outside $A$" then
   \[N^{\Bor_\Gamma^N}\models``\textrm{$a^*$ diagonalizes $M^{\Bor_\Gamma^M}$ outside $A$}".\]
\end{corollary}

Proofs of both Lemmas \ref{EpresMAD} and \ref{randompresMAD} use an argument similar to that of the proof that the respective posets are $\Dbf$-good (the compactness argument for $\Eor$ and $\omega^\omega$-bounding for $\Bor$).

\begin{question}\label{QSuslinDiag}
   Assume $\Sor$ is a Suslin ccc poset coded in $M$ such that
   $M\models$``$\Sor$ is $\Dbf$-good" and $N\models``a^*$ diagonalizes
   $M$ outside $A$". Does one have:

   \[N^{\Sor^N}\models``\textrm{$a^*$ diagonalizes $M^{\Sor^M}$ outside $A$}"?\]
\end{question}

\begin{lemma}[{\cite[Lemma 12]{VFJB11}}]\label{limitMAD}
   Let $\sbf$ be a coherent pair of FS iterations, $A\in V$ a family of infinite subsets of $\omega$ and $\dot{a}^*$ a $\Por_{i_1,0}$-name for an infinite subset of $\omega$ such that
   \[\Vdash_{\Por_{i_1,\xi}}\textrm{``$\dot{a}^*$ diagonalizes $V_{i_0,\xi}$ outside $A$"}\]
   for all $\xi<\pi$. Then, $\Por_{i_0,\pi}\lessdot\Por_{i_1,\pi}$ and $\Vdash_{\Por_{i_1,\pi}}$ ``$\dot{a}^*$ diagonalizes $V_{i_0,\pi}$ outside $A$".
\end{lemma}

The results above are summarized as follows when considering standard 2D-coherent systems.

\begin{theorem}\label{presMAD}

   Let $\matit$ be a standard 2D-coherent system  with $I^\matit=\gamma+1$ an ordinal and $\pi^\matit=\pi$ satisfying (i) and (ii) of Lemma \ref{GenNoNewReals} and, for each $\alpha<\gamma$, let $\dot{a}_\alpha$ be
   a $\Por_{\alpha+1,0}$-name of an infinite subset of $\omega$ such that $\Por_{\alpha+1,0}$ forces that $\dot{a}_\alpha$ diagonalizes $V_{\alpha,0}$ outside $\{\dot{a}_\varepsilon:\varepsilon<\alpha\}$ and $\Por_{\gamma,0}$ forces $\dot{A}=\{\dot{a}_\alpha:\alpha<\gamma\}$ to be an a.d. family. If $\Sor_\xi\in\{\Cor,\Eor\}\cup\ral$ for all $\xi\in S$ then $\Por_{\gamma,\pi}$ forces that $\dot{A}$ is mad and $\afrak\leq|\gamma|$.
\end{theorem}
\begin{proof}
   Lemmas \ref{Nonewreals}, \ref{FixedDiag}, \ref{EpresMAD}, \ref{randompresMAD} and \ref{limitMAD} imply that $\la V_{\alpha,\pi}:\alpha\leq\gamma\ra$ and $A$ satisfy the hypothesis of Corollary \ref{DiagMad}, so $A$ is mad in $V_{\gamma,\pi}$.
\end{proof}

\begin{remark}\label{mathias}
  \begin{enumerate}[(1)]
     \item Other mad families can be considered in this theory of preservation, for instance, the mad family added by an FS iteration of Mathias-Prikry posets. Given an a.d. family $A\subseteq[\omega]^{\aleph_0}$, let $F(A)\subseteq[\omega]^{\aleph_0}$ be the closure of $\{\omega\menos a: a\in A\}\cup\{\omega\menos n:n<\omega\}$ under finite intersections. Note that the generic real $a^*$ added by the Mathias-Prikry poset $\Mor(F(A))$ is almost disjoint from all the members of $A$ and $|a^*\cap x|=\aleph_0$ for every $x\in V\menos\Iwf(A)$. Moreover, $\Mor(F(A))$ forces that $a^*$ diagonalizes $V$ outside $A$. Thus, for an ordinal $\gamma$ with uncountable cofinality, the FS iteration $\langle\Por_\alpha,\Qnm_\alpha\rangle_{\alpha<\gamma}$ with $\Qnm_\alpha=\Mor(F(A\frestr\alpha))$ adds an a.d. family $A=\{a_\alpha:\alpha<\gamma\}$ where each $a_\alpha$ is the Mathias real added by $\Qnm_\alpha$. By Corollary \ref{DiagMad}, $\Por_\gamma$ forces that $A$ is mad.
     \item Any FS iteration of length $\omega_1$ of non-trivial ccc posets adds a mad family of size $\aleph_1$ (so it forces $\afrak=\aleph_1$), actually, the mad family is defined from the Cohen reals added at limit stages. To understand this, it is enough to note that, if $A\in V$ is a countable a.d. family, then $\Cor\simeq\Mor(F(A))$, so any Cohen generic defines an $\Mor(F(A))$-generic.
  \end{enumerate}
\end{remark}

\begin{remark}\label{ufs}
   A version of the previous theorem was originally proved by Brendle and the first author \cite{VFJB11} for a special case where Mathias-Prikry posets with ultrafilters are considered. In the same way, Mathias-Prikry posets can be incorporated into standard 2D-iterations as in Definition \ref{DefStandMatit}. This was done by the third author in \cite{mejia2} to obtain consistency results about the cardinal invariants $\pfrak$, $\sfrak$, $\rfrak$ and $\ufrak$ in relation with those in Cicho\'n's diagram. But thanks to Lemmas \ref{EpresMAD} and \ref{randompresMAD}, and Remark \ref{Remb=a}, the constructions there can be modified to force, additionally, $\bfrak=\afrak$ (like in Theorem \ref{others}).
\end{remark}

The following is a generalization of a result of Steprans \cite{steprans} which shows that the maximal almost disjoint family added by the forcing $\Hor_\kappa$ is indestructible after forcing with some particular posets. Steprans' result can then be deduced when $\kappa=\omega_1$ (so $\Hor_{\omega_1}=\Cor_{\omega_1}$) and $\Qnm_\xi=\Cor$ for all $\xi<\pi$.

\begin{theorem}\label{MainFSit}
   Let $\kappa$ be an uncountable regular cardinal. After forcing with $\Hor_\kappa$, any FS iteration $\la\Por_\xi,\Qnm_\xi\ra_{\xi<\pi}$ where each iterand is either
   \begin{enumerate}[(i)]
      \item in $\{\Cor,\Eor\}\cup\ral$ or
      \item a ccc poset of size $<\kappa$
   \end{enumerate}
   preserves the mad family added by $\Hor_\kappa$.
\end{theorem}
\begin{proof}
We reconstruct the iteration $\Hor_\kappa$ followed by $\la\Por_\xi,\Qnm_\xi\ra_{\xi<\pi}$ as a standard 2D-coherent system $\matit$ so that $\Por^\matit_{\kappa,\xi}=\Hor_\kappa\ast\Por_\xi$ for all $\xi\leq\pi$. The construction goes as follows (see Definition \ref{DefStandMatit}):

\begin{enumerate}
 \item $I^\matit= \kappa+1$ and $\pi^\matit= \pi$.

 \item For each $\alpha\leq\kappa$, $\Por^\matit_{\alpha,0} = \Hor_\alpha$.

 \item The partition $\langle S^\matit, C^\matit \rangle$ of $\pi^\matit$ corresponds to the set of ordinals in the iteration where a poset coming from (i) or (ii) is used. In other words, $\xi \in S^\matit$ if (i) holds for $\Qnm_\xi$, and $\xi \in C^\matit$ otherwise.

 \item The functions $\Delta^\matit: C^\matit \to \kappa$ and the sequences $\langle \Sor^\matit_\xi: \xi \in S^\matit \rangle$ and $\langle \dot{\Qor}^\matit_\xi: \xi \in C^\matit \rangle$ are constructed by recursion on $\xi< \pi$ along with the FS iterations of the 2D-coherent system. We split into the following cases:

 \begin{itemize}
  \item If $\xi \in S^\matit$ define $\Sor^\matit_\xi$ to be one of the posets in the set $\{\Cor,\Eor\}\cup\ral$ depending on what $\Por_{\xi}$ forces $\dot{\Qor}_\xi$ to be.
  \item If $\xi \in C^\matit$ we define both $\Delta^\matit(\xi)$ and $\dot{\Qor}^\matit_\xi$, the latter as a $\Por^\matit_{\Delta^\matit(\xi),0}$-name. Since $\xi \in C^\matit$ we have that $\dot{\Qor}_\xi$ is a $\Por^\matit_{\kappa,\xi}$-name for a ccc poset of size $<\kappa$, hence without loss of generality we can assume that the domain of $\dot{\Qor}_\xi$ is an ordinal $\gamma_\xi< \kappa$ (not just a name). By Lemma \ref{GenNoNewReals}, $\dot{\Qor}_\xi$ is (forced by $\Por^\matit_{\kappa,\xi}$ to be equal to) a $\Por^\matit_{\alpha,\xi}$-name $\dot{\Qor}^\matit_\xi$ for some $\alpha<\kappa$. So put $\Delta^\matit(\xi)=\alpha+1$.\footnote{Though it would be fine to put $\Delta^\matit(\xi)=\alpha$, we prefer $\alpha+1$ because we additionally have that, for any $\gamma<\kappa$ of uncountable cofinality and for any $\xi\leq\pi$, $\R\cap V_{\gamma,\xi}=\R\cap\bigcup_{\alpha<\gamma}V_{\alpha,\xi}$.}
 \end{itemize}
\end{enumerate}
Notice that $\matit$ satisfies the assumptions of Theorem \ref{presMAD} for the mad family $A$ added by $\Hor_\kappa$, so $A$ is still mad in $V^\matit_{\kappa,\pi}$.

\end{proof}

\begin{remark}\label{RemZhang}
   When $\kappa=\omega_1$ in Theorem \ref{MainFSit}, by Remark \ref{mathias}(2) the result still holds when $\Hor_{\omega_1}$ is replaced by any FS iteration of length with cofinality $\omega_1$. This is an alternative (and also a generalization) of Zhang's result \cite{zhang} which states that, under CH, there is a mad family in the ground model which stays mad after an FS iteration of $\Eor$.
\end{remark}

\section{Consistency results on Cicho\'n's diagram}\label{SecMain}

In this section, we prove the consistency of certain constellations in Cicho\'{n}'s diagram where, additionally, the almost disjointness number can be decided (equal to $\bfrak$). For all the results, we fix uncountable regular cardinals $\theta_0\leq\theta_1\leq\kappa\leq\mu\leq\nu$ and a cardinal $\lambda\geq\nu$. We denote the ordinal product between cardinals by, e.g., $\lambda\cdot\mu$.

The following summarizes the results in \cite[Sect. 3]{mejia} but in addition we get that $\bfrak=\afrak$ can be forced.

\begin{theorem}\label{FSappl}
   Assume $\lambda=\lambda^{<\kappa}$ and $\lambda'\geq\lambda$ with $(\lambda')^{\aleph_0}=\lambda'$. For each of the items below, there is a ccc poset forcing the corresponding statement.
   \begin{enumerate}[(a)]
      \item $\add(\Nwf)=\theta_0$, $\cov(\Nwf)=\theta_1$, $\bfrak=\afrak=\non(\Mwf)=\kappa$ and $\cov(\Mwf)=\cfrak=\lambda$.
      \item $\add(\Nwf)=\theta_0$, $\cov(\Nwf)=\theta_1$, $\bfrak=\afrak=\kappa$, $\non(\Mwf)=\cov(\Mwf)=\mu$ and $\dfrak=\non(\Nwf)=\cfrak=\lambda$.
      \item $\add(\Nwf)=\theta_0$, $\bfrak=\afrak=\kappa$, $\cov(\Iwf)=\non(\Iwf)=\mu$ for $\Iwf\in\{\Mwf,\Nwf\}$ and $\dfrak=\cfrak=\lambda$.
      \item $\non(\Nwf)=\aleph_1$, $\bfrak=\afrak=\kappa$, $\dfrak=\lambda$ and $\cov(\Nwf)=\cfrak=\lambda'$.
   \end{enumerate}
\end{theorem}
\begin{proof}
   The proofs are basically the same as in \cite{mejia} combined with the methods of preservation of mad families developed in Section \ref{SecMad}. We sketch these proofs for completeness. For all the items, start adding a mad family with $\Hor_\kappa$.
   \begin{enumerate}[(a)]

   \item Construct an iteration as in the last part of \cite[Thm. 2]{mejia}. To be more precise, perform an FS iteration $\la\Por_\alpha,\Qnm_\alpha\ra_{\alpha<\lambda}$ where each $\Qnm_\alpha$ is either
 \begin{enumerate}[(i)]
      \item a $\sigma$-linked subposet of $\Loc$ of size $<\theta_0$,
      \item a subalgebra of $\Bor$ of size $<\theta_1$ or
      \item a $\sigma$-centered subposet of $\Dor$ of size $<\kappa$.
   \end{enumerate}
        The iteration is constructed by a book-keeping device so that any $\sigma$-linked subposet of $\Loc$ of size $<\theta_0$ that lives in a intermediate step is used in a further step of the iteration. Likewise in relation to (ii) and (iii).

   By Theorem \ref{MainFSit}, $\Por_\lambda$ forces $\afrak\leq\kappa$. On the other hand, by similar arguments as in \cite[Thm. 2]{mejia}, the other equalities are forced. We just show some of them.

   \underline{$\add(\Nwf)=\theta_0$}. The inequality $\add(\mathcal{N}) \leq \theta_0$ follows from both the fact that \\ $\add(\mathcal{N})= \bfrak(\Lc)$ (see Example \ref{ExmpGood}(4)) and that all the posets we are using in the iteration are $\theta_0$-$\Lc$-good, so Theorem \ref{MainPresThm} applies and we get $\bfrak(\Lc) \leq \theta_0$. On the other hand, $\add(\mathcal{N}) \geq \theta_0$ follows from the book-keeping corresponding to (i).



    \underline{$\cov(\Mwf)=\cfrak=\lambda$}. The inequality
    $\cov(\mathcal{M}) \geq \lambda$ is a simple consequence of the
    equality $\cov(\mathcal{M})= \dfrak(\Ed)$ together with Theorem
    \ref{MainPresThm}; on the other hand, $\cfrak\leq\lambda$ because,
    in the ground model, $|\Hor_\kappa\ast\Por_\lambda|\leq\lambda$.

   \item As in (a), perform an FS iteration $\la\Por_\alpha,\Qnm_\alpha\ra_{\alpha<\lambda \cdot \mu}$ as in \cite[Thm. 3]{mejia} where each $\Qnm_\alpha$ is either
   \begin{enumerate}[(i)]
      \item a $\sigma$-linked subposet of $\Loc$ of size $<\theta_0$,
      \item a subalgebra of $\Bor$ of size $<\theta_1$,
      \item a $\sigma$-centered subposet of $\Dor$ of size $<\kappa$ or
      \item $\Eor$.
   \end{enumerate}
   By counting arguments, the FS iteration is constructed so that, for any $\alpha<\mu$, each $\sigma$-linked subposet of $\Loc$ of size $<\theta_0$ living in $V_{\lambda\cdot\alpha}$ is used in the iteration at stage $\lambda\cdot\alpha+\xi$ for some $\xi<\lambda$. Likewise for (ii) and (iii).

   \item Perform an FS iteration $\la\Por_\alpha,\Qnm_\alpha\ra_{\alpha< \lambda \cdot \mu}$ as in \cite[Thm. 4]{mejia}. In this case, each $\Qnm_\alpha$ is either:

   \begin{enumerate}[(i)]
      \item a $\sigma$-linked subposet of $\Loc$ of size $<\theta_0$,
      \item a $\sigma$-centered subposet of $\Dor$ of size $<\kappa$ or
      \item $\Bor$.
   \end{enumerate}

   Counting arguments are used as in (b).

   \item After the iteration in (a) force with $\Bor_{\lambda'}$.
   \end{enumerate}
\end{proof}

Now we turn to prove some consistency results with standard 3D-coherent systems (see Definitions \ref{DefCoherent}(3) and \ref{DefStandMatit}). Recall that, if $\tbf$ is such a system with $I^\tbf=(\gamma+1)\times(\delta+1)$, standard 2D-coherent systems $\tbf_\alpha$ can be extracted for each $\alpha\leq\gamma$ and $\tbf^\beta$ for each $\beta\leq\delta$. When referring to Figure \ref{Figcubeit}, we call the vertical axis the \emph{$\alpha$-axis}, the axis pointing ``perpendicular to the sheet of paper" is the \emph{$\beta$-axis} and the horizontal axis is the \emph{$\xi$-axis.} To get a picture of these 2D-systems, in Figure \ref{Figcubeit}, $\tbf_\alpha$ is the 2D-system obtained by restricting the 3D rectangle to the horizontal plane on $\alpha$ (i.e., fixing $\alpha$ on the $\alpha$-axis), while $\tbf^\beta$ is the restriction to the vertical plane on $\beta$ (i.e., fixing $\beta$ on the $\beta$-axis). These 2D-coherent systems allow us to directly apply the results in the previous sections to 3D-
coherent systems. In consequence, we have the following general result for standard 3D-coherent systems.

\begin{theorem}\label{Maincube}
   Let $\tbf$ be a standard 3D-coherent system with $I^\tbf=(\gamma+1)\times(\delta+1)$ and $\matit$ a standard 2D-coherent system with $I^\matit=\gamma+1$ and $\pi^\matit=\delta$ such that $\Por_{\alpha,\beta,0}=\Por^\tbf_{\alpha,\beta,0}=\Por^\matit_{\alpha,\beta}$ for all $\alpha\leq\gamma$ and $\beta\leq\delta$. Let $\Rbf=\la X,Y,\sqsubset\ra$ be a Polish relational system coded in $V$. Assume
   \begin{enumerate}[(I)]
      \item $\matit$ satisfies the hypotheses of either
           \begin{enumerate}[(i)]
              \item Lemma \ref{GenNoNewReals}(i) and (ii) and Theorem \ref{MainMatit} with $\la \dot{c}_\alpha:\alpha<\gamma\ra$ and $\Rbf$, or
              \item Theorem \ref{presMAD} with $\dot{A}=\{\dot{a}_{\alpha}:\alpha<\gamma\}$
           \end{enumerate}
           (note that, in either case, $\gamma$ has uncountable cofinality),
      \item all the posets that form $\matit$ are non-trivial (see Definition \ref{DefStandMatit}(iii) and (iv)),
      \item all the posets that form $\tbf$ are non-trivial (see Definition \ref{DefStandMatit}(iii) and (iv)),
      \item $\delta$ and $\pi$ have uncountable cofinality,
      \item for $\xi\in S=S^\tbf$, $\Qnm_{\alpha,\beta,\xi}$ is forced to be $\Rbf$-good by $\Por_{\alpha,\beta,\xi}$ for all $\alpha\leq\gamma$ and $\beta\leq\delta$, and
      \item if (I)(ii) is assumed then $\Sor_\xi\in\{\Cor,\Eor\}\cup\mathcal{R}$ for all $\xi\in S$.
   \end{enumerate}
   Then, $\Por_{\gamma,\delta,\pi}$ forces
   \begin{enumerate}[(a)]
      \item $\non(\Mwf)\leq\cf(\pi)\leq\cov(\Mwf)$,
      \item $\bfrak(\Rbf)\leq\min\{\cf(\delta),\cf(\pi)\}\leq\max\{\cf(\delta),\cf(\pi)\}\leq\dfrak(\Rbf)$,
      \item $\bfrak(\Rbf)\leq\min\{\cf(\gamma),\cf(\delta),\cf(\pi)\}\leq
      \max\{\cf(\gamma),\cf(\delta),\cf(\pi)\}\leq\dfrak(\Rbf)$ \underline{when} (I)(i) is assumed and
      \item $\afrak\leq|\gamma|$ \underline{when} (I)(ii) is assumed.
   \end{enumerate}
\end{theorem}
\begin{proof}

   \begin{enumerate}[(a)]
    \item Any FS iteration of length $\pi$ of uncountable cofinality adds cofinally $\cf(\pi)$-many Cohen reals which witness $\non(\mathcal{M})\leq\cf(\pi)\leq\cov(\mathcal{M})$. Also note that the FS iteration $\la\Por_{\gamma, \delta, \xi}, \dot{\Qor}_{\gamma, \delta, \xi}: \xi < \pi\ra$ generates the final extension $V_{\gamma,\delta,\pi}$ of the coherent system $\tbf$.

    \item We look at the 2D-coherent system $\tbf_\gamma$. As the chain of posets $\la \Por_{\gamma,\beta,0}:\beta\leq\delta\ra$ is generated by an FS iteration of ccc posets, for a fixed cofinal sequence $\la\beta_\zeta:\zeta<\cf(\delta)\ra$ in $\delta$ of limit ordinals, for each $\zeta<\cf(\delta)$ there is a $\Por_{\gamma,\beta_{\zeta+1},0}$-name $\dot{c}'_\zeta$ for a Cohen real over $V_{\gamma,\beta_\zeta,0}$. Thus, $\tbf_\gamma$ and $\la\dot{c}'_\zeta:\zeta<\cf(\delta)\ra$ satisfy the hypotheses of Theorem \ref{MainMatit} by (V), so $\Por_{\gamma,\delta,\pi}$ forces $\bfrak(\Rbf)\leq\cf(\delta)\leq\dfrak(\Rbf)$. Besides, since $\bfrak(\Rbf) \leq \non(\mathcal{M})$ and $\cov(\mathcal{M})\leq \dfrak(\Rbf)$, (a) immediately implies $\bfrak(\Rbf)\leq \cf(\pi) \leq\dfrak(\Rbf)$.

    \item We first look at the 2D-coherent system $\matit$. By Theorem \ref{MainMatit}, $\Por_{\alpha+1,\delta,0}$ forces that $\dot{c}_\alpha$ is $\Rbf$-unbounded over $V_{\alpha,\delta,0}$ for every $\alpha<\gamma$. Now, we apply Theorem \ref{MainMatit} to $\tbf^\delta$ to conclude that $\bfrak(\Rbf)\leq\cf(\gamma)\leq\dfrak(\Rbf)$.

    \item By Theorem \ref{presMAD} applied to the 2D-coherent system $\matit$, each $\dot{a}_\alpha$ is forced by $\Por_{\alpha+1,\delta,0}$ to diagonalize $V_{\alpha,\delta,0}$ outside $\dot{A}\frestr\alpha$ for each $\alpha<\gamma$ and furthermore, using the same theorem one more time for the coherent system $\tbf^\delta$, $\Por_{\alpha+1,\delta,\pi}$ forces that $\dot{a}_\alpha$ diagonalizes $V_{\alpha,\delta,\pi}$ outside $\dot{A}\frestr\alpha$. Thus,
        the maximality of $A$ is preserved in $V_{\gamma,\delta,\pi}$ and so $\afrak \leq |\gamma|$.
   \end{enumerate}
\end{proof}

In our applications and in accordance with the previous result, we consider standard 3D-coherent systems where $\la\Por_{\alpha,\beta,0}:\alpha\leq\gamma,\beta\leq\delta\ra$ is generated by a standard 2D-coherent system.

\begin{definition}\label{Defstartmat}
Given ordinals $\gamma$ and $\delta$, define the following standard 2D-coherent systems.
\begin{enumerate}[(1)]

\item The system $\matit^\Cor(\gamma,\delta)$ where

\begin{enumerate}[(i)]
   \item $I^{\matit^\Cor(\gamma,\delta)}=\gamma+1$,
   \item $\Por^{\matit^\Cor(\gamma,\delta)}_{\alpha,0}=\Cor_\alpha$ for each $\alpha\leq\gamma$, and
   \item $\pi^{\matit^\Cor(\gamma,\delta)}=\delta$, $S=\delta$, $C=\emptyset$ and $\Sor_\beta=\Cor$ for all $\beta<\delta$.
\end{enumerate}

\item The system $\matit^*(\gamma,\delta)$ where
\begin{enumerate}[(i)]
   \item $I^{\matit^*(\gamma,\delta)}=\gamma+1$,
   \item $\Por^{\matit^*(\gamma,\delta)}_{\alpha,0}=\Hor_\alpha$ for each $\alpha\leq\gamma$, and
   \item $\pi^{\matit^*(\gamma,\delta)}=\delta$, $S=\delta$, $C=\emptyset$ and $\Sor_\beta=\Cor$ for all $\beta<\delta$.
\end{enumerate}
\end{enumerate}
\end{definition}

If both $\gamma$ and $\delta$ have uncountable cofinality, it is clear that both $\matit^\Cor(\gamma,\delta)$ and $\matit^*(\gamma,\delta)$ satisfy (I) and (II) of Theorem \ref{Maincube}, moreover, the former satisfies (I)(i) and the latter satisfies (I)(ii). These standard 2D-coherent systems are the starting point for the 3D-coherent systems constructed to prove the main results below.

Note that in Theorems \ref{Coflarge}(b), \ref{contlarge}(c) and (d) we cannot say anything about $\afrak$ because full Hechler generics are added (see the discussion about full and restricted generics after Definition \ref{DefStandMatit}) so mad families are not preserved anymore in the way proposed in Section \ref{SecMad}. For these results we start with $\matit^\Cor(\cdot,\cdot)$. For the results where we can force $\bfrak=\afrak$ we start with $\matit^*(\cdot,\cdot)$ (we can start with $\matit^\Cor(\cdot,\cdot)$ as well, but $\afrak$ should be ignored in that case). Observe that the results below are three-dimensional versions of the 2D-coherent systems constructed in \cite[Sect. 6]{mejia}.

\begin{figure}
   \begin{center}
     \includegraphics{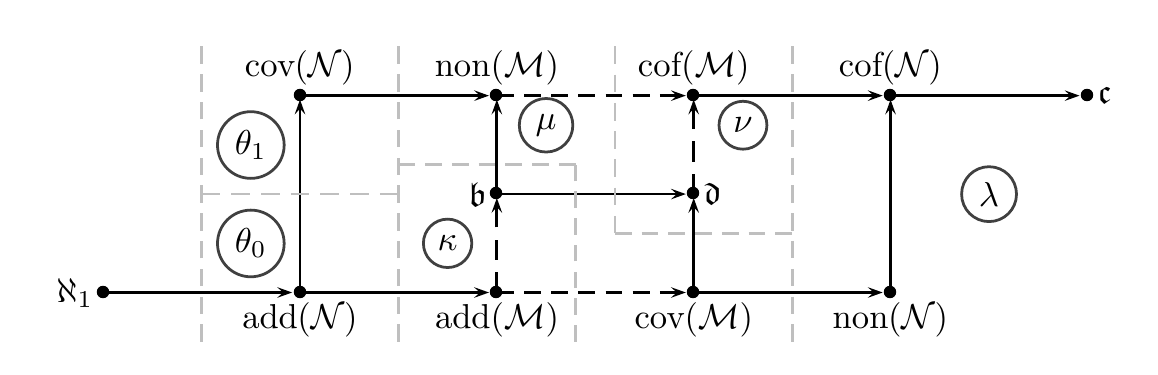}
     \caption{Cicho\'n's diagram as in Theorem \ref{7values}.}
     \label{Fig7values}
   \end{center}
\end{figure}

We first prove that there is a constellation of Cicho\'n's diagram with 7 different values as illustrated in Figure \ref{Fig7values}.

\begin{theorem}\label{7values}
   Assume $\lambda^{<\theta_1}=\lambda$. Then, there is a ccc poset forcing $\add(\Nwf)=\theta_0$, $\cov(\Nwf)=\theta_1$, $\bfrak=\afrak=\kappa$, $\non(\Mwf)=\cov(\Mwf)=\mu$, $\dfrak=\nu$ and $\non(\Nwf)=\cfrak=\lambda$.
\end{theorem}

\begin{proof}

Let $V$ be the ground model where we perform an FS iteration which comes from the standard 3D-coherent system $\tbf$ constructed as follows. Fix a bijection $g=\la g_0,g_1,g_2\ra:\lambda\to\kappa\times\nu\times\lambda$.

\begin{enumerate}[(1)]
 \item $\gamma = \kappa+1$, $\delta = \nu+1$ and $\pi= \lambda \cdot \nu \cdot \mu$.

 \item $\langle \Por_{\alpha, \beta, 0}: \alpha \leq \kappa, \beta\leq \nu \rangle$ is obtained from $\matit^\ast(\kappa, \nu)$.

 \item Consider $\lambda\cdot\nu\cdot\mu$ as the disjoint union of the $\nu \cdot \mu$-many intervals $I_\zeta =[l_\zeta, l_{\zeta + 1})$ (for $\zeta < \nu \cdot \mu$) of order type $\lambda$. Let $S := \{l_\zeta : \zeta < \nu \cdot \mu \}$ and  $C =\pi \smallsetminus S$ (note that $l_\zeta=\lambda\cdot\zeta$).

 \item A function $\Delta=\la\Delta_0,\Delta_1\ra:C \to\kappa\times\nu$ such that the following properties are satisfied:
      \begin{enumerate}[(i)]
        \item For all $\xi <\pi$, both $\Delta_0(\xi)$ and $\Delta_1(\xi)$ are successor ordinals,\footnote{Both ordinals $\Delta_0(\xi)$ and $\Delta_1(\xi)$ are successor because, if they are limits of uncountable cofinality and we force with $\Dor^V_{\Delta(\xi),\xi}$ above $(\Delta(\xi),\xi)$ and trivial otherwise, then $\R\cap V_{\Delta(\xi),\xi+1}$ may not be $\R\cap\bigcup_{\alpha<\Delta_0(\xi),\beta<\Delta_1(\xi)}V_{\alpha,\beta,\xi+1}$.}
        \item $\Delta^{-1}(\alpha+1,\beta+1)\cap\{l_\zeta+1:\zeta<\nu\cdot\mu\}$ is cofinal in $\pi$ for any $(\alpha,\beta)\in\kappa\times\nu$, and
        \item for fixed $\zeta<\nu\cdot\mu$ and $e<2$, $\Delta(l_\zeta+2+2\cdot\varepsilon+e)=(g_0(\varepsilon)+1,g_1(\varepsilon)+1)$ for all $\varepsilon<\lambda$.
      \end{enumerate}
 \item $\Sor_{\xi} = \Eor$ for all $\xi \in S$.

 \item Fix, for each $\alpha<\kappa$, $\beta<\nu$ and $\zeta<\nu\cdot\mu$, two sequences $\langle \dot{\Loc}^\zeta_{\alpha,\beta,\eta} \rangle_{\eta<\lambda}$ and $\langle \dot{\Bor}^\zeta_{\alpha,\beta,\eta} \rangle_{\eta < \lambda}$ of $\Por_{\alpha,\beta, l_\zeta}$-names for all $\sigma$-linked subposets of the localization forcing $\Loc^{V_{\alpha,\beta, l_\zeta}}$ of size $< \theta_0$ and all subalgebras of random forcing $\Bor^{V_{\alpha, \beta, l_\zeta}}$ of size $< \theta_1$, respectively.

     Given $\xi\in C$, define $\Qnm_\xi$ according to the following cases.
     \begin{enumerate}[(i)]
          \item If $\xi= l_\zeta +1$ then $\Qnm_\xi$ is a $\Por_{\Delta(\xi),\xi}$-name for the poset $\Dor^{V_{\Delta(\xi),\xi}}$, the Hechler poset adding a dominating real $\dot{d}_\zeta$ over the model $V_{\Delta(\xi),\xi}$.
          \item If $\xi= l_\zeta+ 2+2\varepsilon$ with $\varepsilon<\lambda$ then $\dot{\Qor}_{\xi} = \dot{\Loc}^\zeta_{g(\varepsilon)}$.
          \item If $\xi= l_\zeta+2+ 2\varepsilon +1$ with $\varepsilon<\lambda$ then $\dot{\Qor}_{\xi} = \dot{\Bor}^\zeta_{g(\varepsilon)}$.
     \end{enumerate}

\end{enumerate}

We prove that $V_{\kappa,\nu,\pi}$ satisfies the statements of this theorem.

\begin{clm}\label{reals}
   If $X\in V_{\kappa,\nu,\pi}$ is a set of reals of size $<\mu$, then
   there are $(\beta,\zeta)\in\nu\times(\nu\cdot\mu)$ so that $X\in
   V_{\kappa,\beta,l_\zeta}$. Furthermore, if $|X|<\kappa$, then there
   is also an $\alpha$ less than $\kappa$ such that $X\in V_{\alpha,\beta,l_\zeta}$.
\end{clm}
\begin{proof}
  As $\cf(\pi)= \mu$ and $V_{\kappa,\nu,\pi}$ is obtained by an FS iteration of length $\pi$, there is a $\zeta<\nu\cdot\mu$ such that $X\in V_{\kappa,\nu,l_\zeta}$ (because $\{l_\zeta:\zeta<\nu\cdot\mu\}$ is cofinal in $\pi$). Now, look at the 2D-coherent system $\tbf_\kappa$ and apply Corollary \ref{Nonewreals} to find a $\beta<\nu$ so that $X\in V_{\kappa,\beta,l_\zeta}$. In the case that $|X|<\kappa$, apply Corollary \ref{Nonewreals} to $\tbf^\beta$ to find an  $\alpha < \kappa$ so that $X$ belongs to $V_{\alpha, \beta, l_\zeta}$.
\end{proof}

\underline{$\add(\mathcal{N}) = \theta_0$}. For the inequality $\add(\mathcal{N}) \geq \theta_0$ take an arbitrary set $X$ of reals in $V_{\kappa,\nu,\pi}$ of size $< \theta_0$ so, by Claim \ref{reals}, there is a triple of ordinals $(\alpha,\beta,\zeta) \in \kappa \times \nu \times (\nu\cdot\mu)$ such that $X\in V_{\alpha,\beta,l_\zeta}$. In $V_{\alpha,\beta,l_\zeta}$, there is a transitive model $N$ of (a large enough finite fragment of) ZFC such that $X\subseteq N$ and $|N|<\theta_0$. Then, there exists an $\eta < \lambda$ such that $\Loc^\zeta_{\alpha,\beta,\eta}=\Loc^N$. Put $\varepsilon= g^{-1}(\alpha, \beta, \eta)$ and $\xi'= l_\zeta +2+2\varepsilon$, so $\Qor_{\xi'}=\Loc^\zeta_{\alpha,\beta,\eta}=\Loc^N$ adds a generic slalom over $N$ and, therefore, it localizes all the reals in $X$.

To obtain the converse inequality, apply Theorem \ref{MainPresThm} to $\la\Por_{\kappa,\nu,\xi},\Qnm_{\kappa,\nu,\xi}\ra_{\xi<\pi}$ and $\theta_0$.

\underline{$\cov(\mathcal{N}) = \theta_1$}. This case is similar to the one above. To get $\cov(\mathcal{N}) \geq \theta_1$ take an arbitrary family $Z$ of Borel null sets coded in $V_{\kappa,\nu,\pi}$ of size $< \theta_1$ so, by Claim \ref{reals}, there exists $(\alpha, \beta, \zeta) \in \kappa \times \nu \times (\nu\cdot\mu)$ such that the sets in $Z$ are already coded in $V_{\alpha, \beta, l_\zeta}$. Hence, as in the previous argument, there exists an ordinal $\eta < \lambda$ such that the generic random real added by $\Bor^\zeta_{\alpha, \beta, \eta }$ avoids all the Borel sets in $Z$. Put $\varepsilon= g^{-1}(\alpha, \beta, \eta)$ and $\xi'= l_\zeta +2+2\varepsilon +1$, so $\Qor_{\xi'}=\Bor^\zeta_{\alpha,\beta,\eta}$ and the random real it adds is already in $V_{\alpha+1, \beta+1, \xi'+1}$.

Conversely, since the posets we use in the FS iteration $\la\Por_{\kappa,\nu,\xi},\Qnm_{\kappa,\nu,\xi}\ra_{\xi<\pi}$ are $\theta_1$-$\Edb$-good posets and $\cov(\mathcal{N}) \leq \bfrak(\Edb)$, Theorem \ref{MainPresThm} implies that, in $V_{\kappa,\nu,\pi}$, $\bfrak(\Edb)\leq \theta_1$.

\underline{$\non(\mathcal{M}) = \cov(\mathcal{M}) = \mu$}. The inequalities $\non(\mathcal{M}) \leq \mu \leq \cov(\mathcal{M})$ follow from Theorem \ref{Maincube}(a). Conversely, from the cofinally $\mu$-many eventually different reals added by the iteration $\la\Por_{\kappa,\nu,\xi},\Qnm_{\kappa,\nu,\xi}\ra_{\xi<\pi}$, we force the inequalities $\cov(\mathcal{M}) \leq \mu$ and $\non(\mathcal{M})\geq \mu$.

\underline{$\add(\mathcal{M}) = \bfrak= \afrak = \kappa$}. Given a family $F$ of reals in $V_{\kappa,\nu,\pi}$ of size $<\kappa$, we can find a $(\alpha, \beta,\zeta) \in \kappa \times \nu \times (\nu\cdot\mu)$ such that $F \in V_{\alpha, \beta,l_\zeta }$. We use now the restricted dominating reals $\{\dot{d}_\zeta : \zeta < \nu\cdot\mu \}$. Since $(\Delta)^{-1}(\alpha+1,\beta+1) \cap \{l_\zeta+1 : \zeta < \nu \cdot \mu\}$ is cofinal in $\pi$, there exists a $\zeta' \in [\zeta,\nu\cdot\mu)$ such that $\Delta(l_{\zeta'}+1)=(\alpha+1,\beta+1)$ and then the real $\dot{d}_{\zeta'}$ added by $\Qor_{\alpha+1,\beta+1,\xi'}$, where $\xi'=l_{\zeta'}+1$, dominates all the reals in $F$.

On the other hand, $\afrak \leq \kappa$ follows from Theorem \ref{Maincube} which guarantees that the mad family added along the $\alpha$-axis, which lives in the model $V_{\kappa,0,0}$, still remains mad in the final extension $V_{\kappa, \nu, \pi}$.

\underline{$\dfrak = \cof(\mathcal{M})= \nu$}. For $V_{\kappa,\nu,\pi} \models \mathfrak{d} \geq \nu$ we just use Theorem \ref{Maincube}. Conversely, to see $V^\Por \models \dfrak \leq \nu$ note that the argument above shows that the family of (restricted) dominating reals $\{\dot{d}_\zeta: \zeta < \nu\cdot\mu\}$ is dominating in $V_{\kappa,\nu,\pi}$.

\underline{$\non(\mathcal{N})= \cof(\mathcal{N}) = \cfrak = \lambda$}. As $\dfrak(\Edb)\leq \non(\mathcal{N})$, from Theorem \ref{MainPresThm} we have that, in $V_{\kappa,\nu,\pi}$, $\dfrak(\Edb)\geq \lvert \pi \lvert = \lambda$. Certainly, $\cfrak\leq\lambda$ holds because $|\Por_{\kappa,\nu,\pi}|=\lambda$.

\end{proof}

\begin{theorem}\label{Coflarge}
   Assume $\lambda^{<\theta_0}=\lambda$. Then, for any of the statements below, there is a ccc poset forcing it.
   \begin{enumerate}[(a)]
       \item $\add(\Nwf)=\theta_0$, $\bfrak=\afrak=\kappa$, $\cov(\Iwf)=\non(\Iwf)=\mu$ for $\Iwf\in\{\Mwf,\Nwf\}$, $\dfrak=\nu$ and $\cof(\Nwf)=\cfrak=\lambda$.
       \item $\add(\Nwf)=\theta_0$, $\cov(\Nwf)=\kappa$, $\add(\Mwf)=\cof(\Mwf)=\mu$, $\non(\Nwf)=\nu$ and $\cof(\Nwf)=\cfrak=\lambda$.
       \item $\add(\Nwf)=\theta_0$, $\cov(\Nwf)=\bfrak=\afrak=\kappa$, $\non(\Mwf)=\cov(\Mwf)=\mu$, $\dfrak=\non(\Nwf)=\nu$ and $\cof(\Nwf)=\cfrak=\lambda$.
    \end{enumerate}
\end{theorem}

\begin{proof}
Fix a bijection $g: \lambda \to \kappa \times \nu \times \lambda$. All the 3D-coherent systems we use in this proof are of the form $\tbf$ where
\begin{enumerate}[(1)]
 \item $\gamma= \kappa+1$, $\delta= \nu+1$ and $\pi= \lambda \cdot \nu \cdot \mu$, the latter of which is the disjoint union of $\nu \cdot \mu$-many intervals $\{I_\zeta:= [l_\zeta, l_{\zeta + 1}) : \zeta < \nu \cdot \mu \}$ of length $\lambda$ where each $l_\zeta:=\lambda\cdot\zeta$.
 \item $S =\{ l_\zeta : \zeta < \nu \cdot \mu\} $ and $C= \pi \smallsetminus S$.
 \item For (a) and (c) $\langle \Por_{\alpha, \beta, 0}: \alpha \leq \kappa, \beta\leq \nu \rangle$ comes from $\matit^\ast(\kappa, \nu)$ and, for (b), it comes from from $\matit^\Cor(\kappa,\nu)$.
 \item A function $\Delta=\la\Delta_0,\Delta_1\ra:C \to\kappa\times\nu$ such that the following properties are satisfied:
      \begin{enumerate}[(i)]
        \item For all $\xi <\pi$, both $\Delta_0(\xi)$ and $\Delta_1(\xi)$ are successor ordinals,
        \item $\Delta^{-1}(\alpha+1,\beta+1)\cap\{l_\eta+1:\eta<\nu\cdot\mu\}$ is cofinal in $\pi$ for each $(\alpha,\beta)\in\kappa\times\nu$; additionally, for (c), $\Delta^{-1}(\alpha+1,\beta+1)\cap\{l_\eta+2:\eta<\nu\cdot\mu\}$ is cofinal in $\pi$ and
        \item for fixed $\zeta<\nu\cdot\mu$, $\Delta(l_\zeta+n_0+\varepsilon)=(g_0(\varepsilon)+1,g_1(\varepsilon)+1)$ for all $\varepsilon<\lambda$, where $n_0=2$ for (a) and (b), and $n_0=3$ for (c).
      \end{enumerate}
 \end{enumerate}

For each of the items below, $\tbf$ is defined appropriately.

\begin{enumerate}[(a)]
  \item For all $\xi\in S$, $\Sor_{\xi}= \Bor$. Fix, for each $\alpha<\kappa$, $\beta<\nu$ and $\zeta<\nu\cdot\mu$, a sequence $\langle \dot{\Loc}^\zeta_{\alpha,\beta,\eta} \rangle_{\eta<\lambda}$ of $\Por_{\alpha,\beta, l_\zeta}$-names for all $\sigma$-linked subposets of $\Loc^{V_{\alpha,\beta, l_\zeta}}$ of size $< \theta_0$. For $\xi \in C$, $\dot{\Qor}_\xi$ is defined according to the following cases.
    \begin{enumerate}[(i)]
       \item If $\xi= l_\zeta +1$ then $\Qnm_\xi$ is a $\Por_{\Delta(\xi),\xi}$-name for the poset $\Dor^{V_{\Delta(\xi),\xi}}$ which adds a dominating real $\dot{d}_\zeta$ over $V_{\Delta(\xi),\xi}$.
       \item If $\xi= l_\zeta+ 2 +\varepsilon$ for some $\varepsilon < \lambda$, then $\dot{\Qor}_{\xi} = \dot{\Loc}^\zeta_{g(\varepsilon)}$.
    \end{enumerate}

  Most of the arguments for each of the cardinal characteristics are identical as the ones presented in Theorem \ref{7values}, so we just present the missing ones.

  \underline{$\non(\mathcal{N})\leq \mu\leq\cov(\Nwf)$}. It holds because we add cofinally $\mu$-many random reals (corresponding to the coordinates $\xi \in S$).

  \underline{$\cof(\mathcal{N})\geq \lambda$}. It is a consequence of both the fact that $\cof(\mathcal{N})= \dfrak(\Lc)$ and Theorem \ref{MainPresThm} which gives us $\dfrak(\Lc) \geq \lvert \pi \lvert= \lambda$.

  \item For all $\xi \in S$, $\Sor_{\xi}= \Dor$ and, for $\xi \in C$, $\dot{\Qor}_\xi$ is defined as in (a) but, in (i), we consider $\Bor^{V_{\Delta(\xi),\xi}}$ instead.

    Recall that, in this construction, our base 2D-coherent system comes from $\matit^\Cor(\kappa,\nu)$. The argument to prove that $V_{\kappa,\nu,\pi}$ satisfies (b) is similar to (a) and to the proof of Theorem \ref{7values}. For instance,

  \underline{$\cov(\mathcal{N})=\kappa$ and $\non(\Nwf)=\nu$}. Given a family $X$ of Borel-null sets coded in $V^\Por$ of size $<\kappa$, we can find $(\alpha, \beta,\zeta) \in \kappa \times \nu \times (\nu\cdot\mu)$ such that all the sets in $X$ are already coded in $V_{\alpha, \beta,l_\zeta }$. Since $\Delta^{-1}(\alpha+1,\beta+1) \cap \{l_\zeta+1 : \zeta < \nu \cdot \mu\}$ is cofinal in $\pi$, there exists $\zeta'\in[\zeta,\lambda)$ such that $\Delta(l_{\zeta'}+1)=(\alpha+1,\beta+1)$ and then the random real $\dot{r}_{\zeta'}$ added by $\dot{\Qor}_{\alpha,\beta,\xi'}$ with $\xi'=l_{\zeta'}+1$ avoids all the sets in $X$. Note that this same argument also proves that the set $\{\dot{r}_{\zeta} : \zeta < \nu\cdot\mu \}$ is not null, so $\non(\mathcal{N}) \leq \nu$.

  Conversely, $\cov(\Nwf)\leq\bfrak(\Edb)\leq\kappa$ and $\nu\leq\dfrak(\Edb)\leq\non(\Nwf)$ are direct consequences of Theorem \ref{Maincube}.

  \underline{$\bfrak=\dfrak=\mu$}. Since the cofinally $\mu$-many dominating reals added by $\la\Por_{\kappa,\nu,\xi},\Qnm_{\kappa,\nu,\xi}\ra_{\xi<\pi}$ forms a scale of length $\mu$.


  \item For all $\xi \in S$, $\Sor_{\xi}= \Eor$. For $\xi \in C$, $\dot{\Qor}_\xi$ is defined according to the following cases
      \begin{enumerate}[(i)]
         \item If $\xi= l_\zeta +1$, then $\Qnm_\xi$ is a $\Por_{\Delta(\xi),\xi}$-name for the poset $\Dor^{V_{\Delta(\xi),\xi}}$.
         \item If $\xi= l_\zeta +2$, then $\Qnm_\xi$ is a $\Por_{\Delta(\xi),\xi}$-name for the poset $\Bor^{V_{\Delta(\xi),\xi}}$.
         \item Otherwise, like (ii) of the proof of (a).
      \end{enumerate}

\end{enumerate}
\end{proof}

\begin{theorem}\label{contlarge}
   Assume $\lambda^{\aleph_0}=\lambda$. Then, for any of the statements below there is a ccc poset forcing it.
   \begin{enumerate}[(a)]
      \item $\add(\Nwf)=\cov(\Nwf)=\bfrak=\afrak=\kappa$, $\non(\Mwf)=\cov(\Mwf)=\mu$, $\dfrak=\non(\Nwf)=\cof(\Nwf)=\nu$ and $\cfrak=\lambda$.
      \item $\add(\Nwf)=\bfrak=\afrak=\kappa$, $\cov(\Iwf)=\non(\Iwf)=\mu$ for $\Iwf\in\{\Mwf,\Nwf\}$, $\dfrak=\cof(\Nwf)=\nu$ and $\cfrak=\lambda$.
      \item $\add(\Nwf)=\cov(\Nwf)=\kappa$, $\add(\Mwf)=\cof(\Mwf)=\mu$, $\non(\Nwf)=\cof(\Nwf)=\nu$ and $\cfrak=\lambda$.
      \item $\add(\Nwf)=\kappa$, $\cov(\Nwf)=\add(\Mwf)=\cof(\Mwf)=\non(\Nwf)=\mu$, $\cof(\Nwf)=\nu$ and $\cfrak=\lambda$.
   \end{enumerate}
\end{theorem}

\begin{proof}
The 3D-coherent systems we use in this proof are of the form $\tbf$ where:

 \begin{enumerate}
 \item $\gamma= \kappa+1$, $\delta= \nu+1$ and $\pi= \lambda \cdot \nu \cdot \mu$ is a disjoint union of $\{I_\zeta= [l_\zeta, l_{\zeta+1}) : \zeta < \nu \cdot \mu \}$ as in Theorem \ref{7values}.
 \item $C=\{l_\zeta:\zeta<\nu\cdot\mu\}$ and $S=\pi\smallsetminus C$.
 \item For items (a) and (b) $\langle \Por_{\alpha, \beta, 0}: \alpha \leq \kappa, \beta\leq \nu \rangle$ comes from $\matit^\ast(\kappa, \nu)$; for (c) and (d), it comes from $\matit^\Cor(\kappa,\nu)$.
 \item A function $\Delta=\la\Delta_0,\Delta_1\ra:C \to\kappa\times\nu$ such that the following properties are satisfied:
      \begin{enumerate}[(i)]
        \item For all $\xi <\pi$, both $\Delta_0(\xi)$ and $\Delta_1(\xi)$ are successor ordinals and
        \item $\Delta^{-1}(\alpha+1,\beta+1)\cap\{l_\zeta:\zeta<\nu\cdot\mu\}$ is cofinal in $\pi$.

      \end{enumerate}
 \end{enumerate}

\begin{enumerate}[(a)]
  \item Put $\Sor_{\xi}= \Eor$ for all $\xi \in S$. For $\xi \in C$, $\dot{\Qor}_\xi=\Loc^{V_{\Delta(\xi),\xi}}$.

  We just prove $\add(\mathcal{N})=\cov(\Nwf)=\bfrak=\kappa$ and $\dfrak=\non(\Nwf)=\cof(\Nwf)=\nu$. If $X$ is a set of reals in $V_{\kappa,\nu,\pi}$ of size $< \kappa$, there is a $(\alpha, \beta, \zeta)\in \kappa\times\nu\times(\nu\cdot\mu)$ such that $X \in V_{\alpha, \beta, l_\zeta}$. Since $\Delta^{-1}(\alpha+1,\beta+1) \cap \{l_\zeta : \zeta < \mu\}$ is  cofinal in $\pi$, there exists a $\zeta'\in[\zeta,\lambda)$ such that $\Delta(l_{\zeta'})=(\alpha+1,\beta+1)$ and then the slalom $\dot{\varphi}_{\zeta'}$ added by $\dot{\Qor}_{\alpha,\beta,l_{\zeta'}}$ localizes all the reals in $X$. Note that $\{\dot{\varphi}_{\zeta}:\zeta<\nu\cdot\mu\}$ witnesses $\cof(\mathcal{N}) \leq \nu$.

  The inequalities $\bfrak,\cov(\Nwf)\leq\kappa$ and $\nu\leq\dfrak,\non(\Nwf)$ follow directly from Theorem \ref{Maincube}.

\item Put $\Sor_{\xi}= \Bor$ for all $\xi \in S$ and, for $\xi \in C$, $\dot{\Qor}_\xi$ is as in (a).

\item Put $\Sor_\xi=\Dor$ for all $\xi \in S$ and, for $\xi \in C$, $\dot{\Qor}_\xi$ is as in (a)

\item For $\xi\in S$, if it is odd then $\Sor_{\xi}= \Dor$, but when it is even then $\Sor_{\xi+1} = \Bor$. For $\xi \in C$, $\dot{\Qor}_\xi$ is defined as in (a).
\end{enumerate}
\end{proof}

We present some other models of constellations of the Cicho\'n diagram known from \cite{mejia}, where additionally $\bfrak=\afrak$ holds.

\begin{theorem}\label{others}
   \begin{enumerate}[(a)]
      \item If $\lambda^{<\theta_1}=\lambda$ then there is a ccc poset forcing $\add(\Nwf)=\theta_0$, $\cov(\Nwf)=\theta_1$, $\bfrak=\afrak=\non(\Mwf)=\kappa$, $\cov(\Mwf)=\dfrak=\nu$ and $\non(\Nwf)=\cfrak=\lambda$.
      \item If $\lambda^{<\theta_0}=\lambda$ then there is a ccc poset forcing $\add(\Nwf)=\theta_0$, $\cov(\Nwf)=\bfrak=\afrak=\non(\Mwf)=\kappa$, $\cov(\Mwf)=\dfrak=\non(\Nwf)=\nu$ and $\cof(\Nwf)=\cfrak=\lambda$.
      \item If $\lambda^{\aleph_0}=\lambda$ then there is a ccc poset forcing $\add(\Nwf)=\non(\Mwf)=\afrak=\kappa$, $\cov(\Mwf)=\cof(\Nwf)=\nu$ and $\cfrak=\lambda$.
   \end{enumerate}
\end{theorem}
\begin{proof}
   For (a) use the construction in \cite[Thm. 20]{mejia}, for (b) see \cite[Thm. 16]{mejia} and for (c) see \cite[Thm. 11]{mejia} but, for the 2D-coherent systems, obtain the first column by forcing with $\Hor_\kappa$ instead.
\end{proof}

\begin{remark}\label{Remb=a}
   By slightly modifying the forcing constructions in Theorems \ref{Coflarge}(b) and \ref{contlarge}(c),(d), it is possible to force, additionally, $\bfrak=\afrak=\mu$. This is thanks to the following idea observed by the anonymous referee, for which we are very grateful. Modify the construction only at steps of the form $\lambda\cdot\nu\cdot\eta$ with $\eta<\mu$. Assume that we have already constructed a $\Por_{0,0,\lambda\cdot\nu\cdot\eta'+1}$-name $\dot{a}_{\eta'}$ of an infinite subset of $\omega$ (this is a Mathias-Prikry generic real added by $\Por_{0,0,\lambda\cdot\nu\cdot\eta'+1}$) for each $\eta'<\eta$, so that $\Por_{0,0,\lambda\cdot\nu\cdot\eta}$ forces that $A\frestr\eta=\{\dot{a}_{\eta'}:\eta'<\eta\}$ is an a.d. family. Put $\Qnm_{\alpha,\beta,\lambda\cdot\nu\cdot\eta}=\Qnm_{0,0,\lambda\cdot\nu\cdot\eta}=\Mor(F(A\frestr\eta))$ for each $\alpha\leq\kappa$, $\beta\leq\nu$ (see Remark \ref{mathias}). Let $\dot{a}_\eta$ be a $\Por_{0,0,\lambda\cdot\nu\cdot\eta+1}$-name of the $\Mor(F(A\frestr\eta))$-
generic real. Note that this real is also $\Mor(F(A\frestr\eta))$-generic over $V_{\kappa,\nu,\lambda\cdot\nu\cdot\eta+1}$ because $F(A\frestr\eta)$ does not depend on $\alpha$ and $\beta$, and the generic real with respect to any $V_{\alpha,\beta,\lambda\cdot\nu\cdot\eta}$ is essentially the same. Thus, as in Remark \ref{mathias}, $\Por_{\kappa,\lambda,\lambda\cdot\nu\cdot\mu}$ forces that $A=\{a_\eta:\eta<\mu\}$ is a mad family. On the other hand, as $\Mor(F(A\frestr\eta))$ is $\sigma$-centered, the arguments to deduce the values of the other cardinal invariants remain intact.
\end{remark}

\begin{remark}\label{RemMA}
   In Theorems \ref{FSappl}, \ref{7values}, \ref{Coflarge} and \ref{others}(a) and (b) we can slightly modify the constructions to force, additionally, $\mathrm{MA}_{<\theta_0}$. For instance, in (6) of the proof of Theorem \ref{7values} we use, instead of $\la\Loc^\zeta_{\alpha,\beta,\eta}\ra_{\eta<\lambda}$, an enumeration $\la\Qnm^\zeta_{\alpha,\beta,\eta}\ra_{\eta<\lambda}$ of all the (nice) $\Por_{\alpha,\beta,l_\zeta}$-names for all the ccc posets with domain an ordinal $<\theta_0$. In (6)(ii), $\Qnm_\xi=\Qnm^\zeta_{g(\epsilon)}$ whenever $\Por_{\kappa,\nu,\xi}$ forces $\Qnm^\zeta_{g(\epsilon)}$ to be ccc, otherwise, $\Qnm_\xi$ is just a name for the trivial poset. In a similar way, we can additionally force $\mathrm{MA}_{<\kappa}$ in Theorems \ref{contlarge} and \ref{others}(c).
\end{remark}


\section{$\Delta^1_3$ well-orders of the reals}\label{SecDelta13}

There has been significant interest in the study of possible
constellations among the classical cardinal characteristics of the
continuum
in the presence of a projective, in fact $\Delta^1_3$-definable,
well-order of the reals (see~\cite{VFSF,VFSFLZ11,VFSFYK14}). Answering
a question of~\cite{VFSFYK14}, we show that each of the constellations
described in the previous section is consistent with the existence of
such a projective
well-order. Since the proofs for the different constellations are very similar, we will only outline
the proof of the following theorem.

\begin{theorem}\label{7values_wo}
In $L$, let $\theta_0<\theta_1<\kappa<\mu<\nu<\lambda$ be uncountable regular cardinals and, in addition, $\lambda <\aleph_\omega$. Then there is a cardinal preserving forcing extension of the constructible universe, $L$, in which there is a $\Delta^1_3$ well-order of the reals and in addition $\add(\Nwf)=\theta_0$, $\cov(\Nwf)=\theta_1$, $\bfrak=\afrak=\kappa$, $\non(\Mwf)=\cov(\Mwf)=\mu$, $\dfrak=\nu$ and $\non(\Nwf)=\cfrak=\lambda$.
\end{theorem}

For convenience we fix natural numbers $n_1<\cdots < n_6$ such that $\theta_0=\omega_{n_1}$, $\theta_1=\omega_{n_2}$, $\kappa=\omega_{n_3}$, $\mu=\omega_{n_4}$, $\nu=\omega_{n_5}$, $\lambda=\omega_{n_6}$. We will work over the constructible universe $L$ and we will use the method of almost disjoint coding as it is developed in~\cite{VFSFLZ11}. We will use the following two notions. We will say that a transitive $\mathrm{ZF}^-$ model $\mathcal{M}$ is {\emph{suitable}} if $\omega_{n_6}^\M$ exists and $\omega_{n_6}^\M=\omega_{n_6}^{L^\M}$ (here $\mathrm{ZF}^-$ denotes $\mathrm{ZF}^-$ minus the power set axiom). For subsets $x,y$ of $\omega$, let $x*y=\{2n:n\in x\}\cup\{2n+1: n\in y\}$ and let $\Box(x)=\{2n+2: n\in x\}\cup\{2n+1: n\notin x\}$. Note that if $x,y$ and $a,b$ are pairs of subsets of $\omega$ such that $\Box(x*y)\subseteq\Box(a*b)$, then $x=a$ and $y=b$.

We will start with a general outline of the proof of
Theorem~\ref{7values_wo}. Our forcing construction can be viewed as a
two stage process: a preliminary stage in which we prepare the
universe, followed by a coding stage in which we will not only adjoin
a well-order of the reals with a $\Delta^1_3$-definition, but also
provide the desired constellations of the cardinal
characteristics. The second stage of our forcing construction
recursively adjoins a well-order of the reals, which we denote $``<"$
and for which we will give an explicit definition later. To give a
$\Delta^1_3$-definition of this well-order, we will make use of a
nicely definable sequence $\bar{S}=\la S_\alpha:\alpha<\pi\ra$
(where $\pi=\lambda\cdot\nu\cdot\mu$) of
stationary, co-stationary subsets of $\omega_{n_6-1}$. Once we fix
such a sequence, for each $\alpha<\pi$ we will adjoin a closed
unbounded subset $C_\alpha$ of $\omega_{n_6-1}$ which is
disjoint from $S_\alpha$. Then with the help of $n_6-2$ many almost
disjoint codings,  we will encode each club $C_\alpha$ into a subset
$X_\alpha$ of $\omega_1$.\footnote{The sets $X_\alpha$ will encode
  also additional information, which is necessary for our
  construction.} Finally we will guarantee that the above
kill of stationarity
is accessible to all countable suitable models: using localizing
posets, we will adjoin the characteristic functions of subsets
$Y_\alpha$ of $\omega_1$, such that $Y_\alpha$ codes $X_\alpha$ and
such that every countable suitable model containing an initial segment
of $Y_\alpha$ will encode an appropriate ``local version'' of a kill
of stationarity. With this, the
first stage of our construction will be complete and we will
denote by $V$ the corresponding generic extension of $L$. The coding
stage will be a modification of the construction providing
Theorem~\ref{7values}, a modification which will allow us to adjoin the
desired $\Delta^1_3$-definition. For every pair of
reals $x,y$ such that $x<y$,  we will generically adjoin a real $r$, which almost disjointly codes the sets
$Y_{\alpha+m}$ for $m\in \Box(x*y)$ (here $\alpha$ will be given by a bookkeeping
function). Thus in particular, $r$ will code the inequality $x<y$ by
encoding a pattern of stationarity, non-stationarity  for the sequence
$\la S_{\alpha+m}:m\in\omega\ra$. A key feature of the entire
construction is the fact that the coherent system of iterations which
we use to provide the final generic extension does not add reals which
accidentally encode a kill of stationarity. This leads to the following
$\Delta^1_3$-definition of the well-order: $x<y$ if and only if there
is a real $r$ such that for every countable suitable model $\mathcal{M}$
containing $r$ there is an ordinal $\alpha<\pi^\mathcal{M}$ such that
$(L[r])^\mathcal{M}\vDash \forall m\in\omega( S_{\alpha+m}^\mathcal{M}\;\hbox{is
  non-stationary iff }m\in\Box(x*y))$.

Now we turn to a more detailed account of the construction. Let $\pi=\lambda\cdot\nu\cdot\mu$ and let $f:\pi\to\lambda$ be a canonical bijection.
For each $\alpha<\pi$, let $W_\alpha$ be the $L$-least subset of $\omega_{n_6-1}$ coding $f(\alpha)$. In the following, we will refer to $W_\alpha$ as the {\emph{$L$-least code of $\alpha$ modulo $f$}}, or simply the {\emph{$L$-least code of $\alpha$}}. We start with a nicely definable sequence $\bar{S}=\la S_\alpha:\alpha<\pi\ra$ of stationary, co-stationary subsets of $\omega_{n_6-1}$. Using bounded approximations, for each $\alpha<\pi$, we add a closed unbounded subset $C_\alpha$ of $\omega_{n_6-1}$ which is disjoint from $S_\alpha$.\footnote{For each $\alpha<\pi$ take $\PP^0_\alpha$ to be the poset of all bounded subsets of $\omega_{{n_6}-1}$ with extension relation end-extension and then take $\PP^0=\prod_{\alpha<\pi}\PP^0_\alpha$ with supports of size $<\omega_{{n_6}-1}$.} Following the notation of~\cite{VFSFLZ11}, for a set of ordinals $X$, $\mathrm{Even}(X)$
denotes the set of all even ordinals in $X$. Now, reproducing the ideas of~\cite{VFSFLZ11}, we can find subsets $Z_\alpha\subseteq\omega_{n_6-1}$ such that

\bigskip
\noindent
$(*)_\alpha$: If $\beta<\omega_{n_6-1}$ and $\mathcal{M}$ is a
suitable model such that $\omega_{n_6-2}\subseteq \mathcal{M}$,
$\omega_{n_6-1}^{\M}=\beta$,
$Z_\alpha\cap\beta\in\mathcal{M}$, then $\mathcal{M}\vDash
\psi(\omega_{n_6-1}^{\M}, Z_\alpha\cap\beta)$, where
$\psi(\omega_{n_6-1}^{\M},X)$ is the formula ``$\mathrm{Even}(X)$ codes a
triple $(\bar{C},\bar{W},\bar{\bar{W}})$ where $\bar{W}$,
$\bar{\bar{W}}$ are the $L$-least codes modulo $f^{\M}$ of ordinals
$\bar{\alpha}$, $\bar{\bar{\alpha}}<\pi^{\M}=
\omega_{n_6}^{\M}\cdot
\omega_{n_5}^{\M}\cdot
\omega_{n_4}^{\M}$ respectively
such that $\bar{\bar{\alpha}}$ is the largest limit ordinal not
exceeding $\bar{\alpha}$, and $\bar{C}$ is a club in $\omega_{n_6-1}^{M}$
disjoint from $S_{\bar{\alpha}}^{\M}$".

\bigskip
For each $m=1,\cdots,n_6-2$, let $\bar{S}^m=\la S^m_\alpha:\alpha<\omega_{n_6-m}\ra$ be a nicely definable in $L_{\omega_{n_6}-m-1}$ sequence of almost disjoint subsets of $\omega_{n_6-m-1}$. Successively using almost disjoint coding with respect to the sequences $\bar{S}^m$ (see~\cite{VFSFLZ11}), we can code the sets $Z_\alpha$ into subsets $X_\alpha$ of $\omega_1$ with the following property.\footnote{Take $X^1_\alpha:=Z_\alpha$ and for each $m=1,\cdots,n_6-2$, let $\PP^m_\alpha$ be the almost disjoint coding of $X^m_\alpha$ via $\bar{S}^m$ (into $X^{m+1}_\alpha$). Define $\PP^{1,m}=\prod_{\alpha<\pi}\PP^m_\alpha$ with supports of size $<\omega_{n_6-m-1}$ and take $\PP^1=\PP^{1,1}*\cdots*\PP^{1,n_6-2}$. Note that $X_\alpha=X_\alpha^{n_6-1}$, the generic of $\PP^{n_6-2}_\alpha$ for each $\alpha$.}

\bigskip
\noindent
$(**)_\alpha$: If $\omega_1<\beta\leq\omega_2$ and $\M$ is a suitable model with $\omega_2^\M=\beta$, $\{X_\alpha\}\cup\omega_1\subseteq\M$, then
$\mathcal{M}\vDash \varphi(\omega_{n_6-1}^{\M}, X_\alpha)$, where
$\varphi(\omega_{n_6-1}^{\M},X)$ is the formula: ``Using the sequences
$(\{\bar{S}^m\}_{m=1}^{m=n_6-2})^{\M}$, the set $X$ almost disjointly codes a
subset $Z$ of $\omega_{n_6-1}^{\M}$ whose even part codes the triple
$(\bar{C},\bar{W},\bar{\bar{W}})$ where $\bar{W}$, $\bar{\bar{W}}$ are
the $L$-least codes modulo $f^{\M}$ of ordinals $\bar{\alpha}$,
$\bar{\bar{\alpha}}<\pi^{\M}$, respectively, such that $\bar{\bar{\alpha}}$
is the largest limit ordinal not exceeding $\bar{\alpha}$, and
$\bar{C}$ is a club in $\omega_{n_6-1}^{\M}$ disjoint from
$S_{\bar{\alpha}}^{\M}$".

\bigskip

Finally, using the posets $\mathcal{L}(X_{\alpha+m}, X_\alpha)$ for $\alpha\in\Lim(\pi)$ (for a set of ordinals $C$, $\Lim(C)$ denotes the set of limit ordinals in $C$), $m\in\omega$ from \cite[Definition 1]{VFSFLZ11}, we can add the characteristic functions of subsets $Y_{\alpha+m}$ of $\omega_1$ such that:\footnote{Take $\PP^2=\prod_{\alpha\in\hbox{Lim}(\pi)}\prod_{m\in\omega}\mathcal{L}(X_{\alpha+m},X_\alpha)$ with countable supports.}

\bigskip
\noindent
$(***)_{\alpha+m}$: If $\beta<\omega_1$, $\mathcal{M}$ is suitable with $\omega_1^\Mwf=\beta$, $Y_{\alpha+m}\cap\beta\in\M$, then
$\mathcal{M}\vDash \varphi(\omega_{n_6-1}^{\M},
X_{\alpha+m}\cap\beta)\wedge \varphi(\omega_{n_6-1}^{\M},
X_\alpha\cap\beta)$.

\bigskip
With this, the preliminary stage of the construction is complete.
We denote by $\PP_0$ the finite iteration of forcing notions described above,
that is $\PP_0=\PP^0*\PP^1*\PP^2$. Note that $\PP_0$ is $\omega$-distributive (the proof is almost identical to~\cite[Lemma 1]{VFSFLZ11}) and so in particular $\PP_0$ does not add new reals. Let $V=L^{\PP_0}$ and let $\bar{\mathcal{B}}=\la B_{\zeta,m}: \zeta<\omega_1, m\in\omega\ra \in L$ be a nicely definable sequence of almost disjoint subsets of $\omega$. As in the proof of Theorem~\ref{7values} partition $\pi$ into intervals $I_\zeta=[l_\zeta,l_{\zeta+1})$ for $\zeta<\nu\cdot\mu$, where $l_\zeta=\lambda\cdot\zeta$, and let
$$C_0=\{2\cdot\zeta'+1:\zeta'<\nu\cdot\mu\}.$$
Furthermore let $C_0^*=\bigcup\{[l_\zeta,l_{\zeta+1}): \zeta\in C_0\}$, let $S^*=\{l_\zeta:\zeta\in\nu\cdot\mu\smallsetminus C_0\}$ and let
$C^*_1=\pi\smallsetminus (S^*\cup \hbox{Lim}(C_0^*))$.

Modifying the 3D-coherent system from the proof of Theorem~\ref{7values}, we will define in $V=L^{\PP_0}$ a standard 3D-coherent system $\tbf^*$ where $\gamma^{\tbf^*}=\kappa+1$, $\delta^{\tbf^*}=\nu+1$, $\pi^{\tbf^*}=\pi$, $S^{\tbf^*}:=S^*$, $C^{\tbf^*}=C^*=\pi\smallsetminus S^*$. The sole difference between $\tbf$ of Theorem \ref{7values} and $\tbf^*$ is the $\xi$-th step of the FS iterations of which $\tbf^*$ consists, when $\xi\in\hbox{Lim}(C_0^*)$. For notational simplicity, $\Por^*_{\alpha,\beta,\xi}=\Por^{\tbf^*}_{\alpha,\beta,\xi}$, $\Qnm^*_{\alpha,\beta,\xi}=\Qnm^{\tbf^*}_{\alpha,\beta,\xi}$, $V^*_{\alpha,\beta,\xi}=V^{\tbf^*}_{\alpha,\beta,\xi}$, $\Delta^*=\Delta^{\tbf^*}:\pi\to\kappa\times\nu$ and so on, while without the asterisk we refer to the components of $\tbf$, that is, $\Por_{\alpha,\beta,\xi}=\Por^{\tbf}_{\alpha,\beta,\xi}$ and so on. Note that $\hbox{Lim}(C_0^*)\subseteq C^*$ and so in particular for $\xi\in\hbox{Lim}(C_0^*)$ we will be adjoining restricted generic reals.

The starting point at $\xi=0$ for $\tbf^*$ is the same as for $\tbf$, that is, $\Por^*_{\alpha,\beta,0}=\Por_{\alpha,\beta,0}$ for all $\alpha\leq\kappa$ and $\beta\leq\nu$. The tasks achieved by the posets $\Qnm_{\alpha,\beta,\xi}$ for $\xi\in S$ (in the notation of the proof of Theorem~\ref{7values}) can be achieved by the corresponding posets $\Qnm^*_{\alpha,\beta,\xi}$ in our modified construction for $\xi\in S^*$, and similarly the tasks achieved by the posets $\Qnm_{\alpha,\beta,\xi}$ for $\xi\in C$ can be accomplished by the posets $\Qnm^*_{\alpha,\beta,\xi}$ for $\xi\in C_1^*$.
Thus, in order to complete the proof of Theorem~\ref{7values_wo}, we are left with  describing the $\xi$-th step for $\xi\in\hbox{Lim}(C^*_0)$ of this modified
construction. It is useful to think of $\sigma^*_\xi=\la \PP_{i,\xi}:i\in I^{\bf{t}^*}\ra$, for $\xi\in\hbox{Lim}(C^*_0)$ as a {\emph{coding section}} of the 3D-coherent system. The reason is that the iterands $\la \QQ^*_{\Delta^*(\xi),\xi}:\xi\in\hbox{Lim}(C^*_0)\ra$ (and correspondingly $\QQ^*_{i,\xi}$ for $i\geq\Delta^*(\xi)$) will be used to introduce a $\Delta^1_3$-definition of a well-order of the reals, which is to be recursively defined along the iteration.
First, we describe this natural well-order of the reals, which arises not only in the modified construction which we are to define, but also in every coherent system we have considered so far in this paper, provided that the corresponding forcing construction is done over the constructible universe $L$.

Our modified 3D-iteration will have the property that for $\alpha^*\leq\kappa$, $\beta^*\leq\nu$ and $\xi^*\leq\pi$, if $G_{\alpha^*,\beta^*,\xi^*}$ is a $\PP_0\ast\PP^*_{\alpha^*,\beta^*,\xi^*}$-generic filter over $L$ then
\begin{multline*}
  L[G_{\alpha^*,\beta^*,\xi^*}] \cap \R = L[ \{\dot{a}_\alpha[G_{\alpha+1,0,0}] :\alpha<\alpha^* \} \cup \{\dot{c}_\beta[G_{0,\beta+1,0}] : \beta<\beta^* \}\\
  \cup \{\dot{u}_{\alpha^*,\beta^*,\xi}[G_{\alpha^*,\beta^*,\xi+1}] : \xi<\xi^* \} ] \cap \R,
\end{multline*}
where $\{\dot{a}_\alpha:\alpha<\kappa\}$ is (the set of names of) the mad family added by $\Hor_\kappa$, $\dot{c}_\beta$ is the Cohen real added by $\Por_{\alpha,\beta+1,0}$ (which does not depend on $\alpha$) and
$\dot{u}_{\alpha,\beta,\xi}$ is a $\PP_0\ast\PP^*_{\alpha,\beta,\xi+1}$-name for the generic real added by $\dot{\Qor}_{\alpha,\beta,\xi}$. Note that, for $\xi\in S^*$, $\PP_0\ast\PP^*_{\alpha,\beta,\xi+1}$ forces $\dot{u}_{\alpha,\beta,\xi}=\dot{u}_{0,0,\xi}$ and, for $\xi\in C^*$, if $\alpha\geq\Delta^*_0(\xi)$ and $\beta\geq\Delta^*_1(\xi)$ then $\PP_0\ast\PP^*_{\alpha,\beta,\xi+1}$ forces $\dot{u}_{\alpha,\beta,\xi}=\dot{u}_{\Delta^*(\xi),\xi}$, otherwise, $\dot{u}_{\alpha,\beta,\xi}$ is just forced to be $\emptyset$. Thus, we only need to look at $\dot{u}_\xi:=\dot{u}_{0,0,\xi}$ when $\xi\in S^*$ and to $\dot{u}_\xi:=\dot{u}_{\Delta^*(\xi),\xi}$ when $\xi\in C^*$.

By recursion on $\alpha^*\leq\kappa$, $\Por_0\ast\Por^*_{\alpha^*,0,0}$ forces that there is a well-order of the reals $\dot{<}_{\alpha^*,0,0}$ which depends only on $\{\dot{a}_\alpha:\alpha<\alpha^*\}$ such that it has $\dot{<}_{\alpha,0,0}$ as an initial segment for every $\alpha<\alpha^*$; by recursion on $\beta^*\leq\nu$, for every $\alpha^*\leq\kappa$, $\Por_0\ast\Por^*_{\alpha^*,\beta^*,0}$ forces that there is a well-order of the reals $\dot{<}_{\alpha^*,\beta^*,0}$ which depends only on $\{\dot{a}_\alpha:\alpha<\alpha^*\}\cup\{\dot{c}_\beta:\beta<\beta^*\}$ such that it has $\dot{<}_{\alpha^*,\beta,0}$ as an initial segment for every $\beta<\beta^*$ and it contains $\dot{<}_{\alpha,\beta^*,0}$ (not necessarily as an initial segment) for every $\alpha<\alpha^*$; and by recursion on $\xi^*\leq\pi$, for all $\alpha^*\leq\kappa$ and $\beta^*\leq\nu$, $ \Por_0\ast\Por^*_{\alpha^*,\beta^*,\xi^*}$ forces that there is a well-order of the reals $\dot{<}_{\alpha^*,\beta^*,\xi^*}$ depending only on $\{\dot{a}_\alpha:
\alpha<\alpha^*\}\cup\{\dot{c}_\beta:\beta<\beta^*\}\cup\{\dot{u}_{\alpha^*,\beta^*,\xi} : \xi<\xi^* \}$ so that it has $\dot{<}_{\alpha^*,\beta^*,\xi}$ as an initial segment for all $\xi<\xi^*$ and contains $\dot{<}_{\alpha,\beta,\xi^*}$ (not necessarily as an initial segment) for every $\alpha\leq\alpha^*$ and $\beta\leq\beta^*$. We denote $\dot{<}_{\xi^*}=\dot{<}_{\kappa,\nu,\xi^*}$. Therefore,  $\PP_0\ast\PP^*_{\kappa,\nu,\xi^*}$ forces that $\dot{<}_{\xi}$ is an initial segment of $\dot{<}_{\xi^*}$ for all $\xi<\xi^*$ and $\PP_0\ast\PP^*_{\kappa,\nu,\pi}$ forces
$$\dot{<}_\pi=\bigcup\{\dot{<}_\xi:\xi<\pi\}$$
which will be the name of the desired well-order. Our modified construction will be done in such a way, that in $L^{\PP_0*\PP^*_{\kappa,\nu,\pi}}$
the reals $\la \dot{u}_\xi[G]:\xi\in\hbox{Lim}(C^*_0)\ra$ will give rise to a $\Delta^1_3$-definition for the well-order $\dot{<}_\pi[G]$.

Now, we turn to the precise definition of the iterands $\QQ^*_\xi$ for $\xi\in \hbox{Lim}(C^*_0)$. We will work in $V$. For each $\xi\in \hbox{Lim}(\pi)$, we will define a $\PP_{\kappa,\nu,\xi}$ name $\dot{A}_\xi$ for a subset of $[\xi,\xi+\omega)$.  Similarly to the construction in~\cite{VFSFLZ11}, for each $\epsilon\in[\omega_{n_6},\omega_{n_6+1})$, fix (in $L$) a bijection $i_\epsilon:\{\la \xi_0,\xi_1\ra: \xi_0<\xi_1<\epsilon\}\to\hbox{Lim}(\omega_{n_6})$.
Fix $\xi\in \hbox{Lim}(C^*_0)$. Then $\xi=l_\zeta+\eta$ for some $\zeta\in C_0$ and $\eta<\omega_{n_6}(=\lambda)$. Suppose $\PP^*_{\alpha,\beta,\xi}$ has been defined for all $\alpha\leq\kappa$, $\beta\leq\nu$. Consider the $\PP^*_{\kappa,\nu,l_\zeta}$-names $\dot{\xi_0},\dot{\xi_1}$ of ordinals for which it is forced that $\la \dot{\xi}_0,\dot{\xi}_1\ra= i_{\hbox{o.t.}(\dot{<}_{l_\zeta})}^{-1}(\eta)$. Furthermore, let $\dot{A}_\xi$ be the $\Por^*_{\kappa,\nu,l_\zeta}$-name of $\xi+(\omega\smallsetminus\Box(x^\zeta_{\dot{\xi}_0}*x^\zeta_{\dot{\xi}_1}))$, where $x^\zeta_\rho$ is the $\rho$-th real in $L[G_{\kappa,\nu,l_\zeta}]\cap [\omega]^{\aleph_0}$ according to the well-order
$\dot{<}_{l_\zeta}$.  By Corollary \ref{Nonewreals}, there are $\alpha<\kappa$ and $\beta<\nu$ such that $\dot{\xi}_0$, $\dot{\xi}_1$ and $\dot{A}_\xi$ are $\PP^*_{\alpha,\beta,l_\zeta}$-names.
Put $\Delta^*(\xi)=(\alpha+1,\beta+1)$ and
$$\dot{\QQ}^*_{\xi}:=\Bigg\{\la s_0,s_1\ra\in[\omega]^{<\aleph_0}\times\Bigg[\bigcup_{m\in\Box(x^\zeta_{\dot{\xi}_0}*x^\zeta_{\dot{\xi}_1})} Y_{\xi+m}\times\{m\}\Bigg]^{<\aleph_0}\Bigg\}$$

\noindent
where $\la t_0,t_1\ra\leq \la s_0, s_1\ra$ if and only if $s_1\subseteq t_1$, $s_0$ is an initial segment of $t_0$ and $(t_0\smallsetminus s_0)\cap B_{\chi,m}=\emptyset$ for all $\la\chi, m\ra\in s_1$.

Note that the real $u_\xi=u_{\Delta^*(\xi),\xi}$ adjoined by $\QQ^*_\xi$ almost disjointly via the sequence $\bar{\mathcal{B}}$ codes the sets $Y_{\xi+m}$ for $m\in \Box(x^\zeta_{\xi_0}*x^\zeta_{\xi_1})$. That is, for every $m\in\Box(x^\zeta_{\xi_0}*x^\zeta_{\xi_1})$, we have $\chi\in Y_{\xi+m}$ iff $|u_\xi\cap B_{\chi,m}|<\omega$. Consider a countable suitable model $\mathcal{M}$ containing $u_\xi$ and let $\mathcal{N}:= (L[u_\xi])^\mathcal{M}$. Then
$\mathcal{N}$ is a suitable countable model, $\omega_1^\mathcal{M}=\omega_1^\mathcal{N}$ and furthermore
$Y_{\xi+m}\cap \omega_1^{\mathcal{N}}\in\mathcal{N}$ for each $m\in\Box(x^\zeta_{\xi_0}*x^\zeta_{\xi_1})$ and so by $(***)_{\xi+m}$ we have
$$\mathcal{N}\vDash \varphi(\omega_{n_6-1}^{\mathcal{N}},
X_{\xi+m}\cap\omega_1^{\mathcal{N}})\wedge \varphi(\omega_{n_6-1}^{\mathcal{N}},
X_\xi\cap\omega_1^{\mathcal{N}}).$$
Then by $(**)_{\xi+m}$ and $(**)_\xi$, we can conclude that there is a $\bar{\bar{\xi}}<\pi^\mathcal{N}$ such that for each $m\in\Box(x^\zeta_{\xi_0}*x^\zeta_{\xi_1})$, $S^\mathcal{N}_{\bar{\bar{\xi}}+m}$ is non-stationary.
To describe the above property of the real $u_\xi$, we will say that {\emph{$u_\xi$ codes a stationarity pattern for $\Box(x^\zeta_{\xi_0}*x^\zeta_{\xi_1})$}}.\footnote{The properties of $\PP_0*\PP^*_{\kappa,\nu,\pi}$ will guarantee that for
$m\in\omega\setminus\Box(x^\zeta_{\xi_0}*x^\zeta_{\xi_1})$, $S^\mathcal{N}_{\bar{\bar{\xi}}+m}$ is stationary in $\mathcal{N}$. Thus the stationarity, non-stationarity pattern of $\la S_{\bar{\bar{\xi}}+m}: m\in\omega \ra$ in $\mathcal{N}$ exactly codes the ordered pair $\la x^\zeta_{\xi_0}, x^\zeta_{\xi_1}\ra$.}

This completes the construction of the modified standard 3D-coherent system. In addition, for every $\xi\in\hbox{Lim}(\pi)\backslash \hbox{Lim}(C_0^*)$ define $\dot{A}_\xi$ to be the canonical $\PP^*_{0,0,\xi}$-name for the interval $[\xi,\xi+\omega)$. Since all posets used to control the cardinal characteristics
in Theorem~\ref{7values} are $\sigma$-linked, as well as the one we used in our coding sections, one can reproduce the proof of~\cite[Lemma 3]{VFSFLZ11} to show that if $G$ is $\PP_0*\PP^*_{\kappa,\nu,\pi}$-generic over $L$, then for each $\eta\in\bigcup_{\xi\in\hbox{Lim}(\pi)} \dot{A}_\xi[G]$ there is no real in $L[G]$
encoding a closed unbounded set disjoint from $S_\eta$. For brevity, we will say that in our final generic extension {\emph{there is no accidental coding of a kill of stationarity by a real}}. This leads to the following $\Sigma^1_3$-definition of $<_\pi$. Let $G$ be a $\PP_0*\PP^*_{\kappa,\nu,\pi}$-generic over $L$ and let $x,y$ be reals in $L[G]$. Then:

\bigskip
\noindent
$x\dot{<}_\pi[G] y$ iff there is a real $r$ such that for every countable suitable model $\mathcal{M}$ such that $r\in\mathcal{M}$, there is $\bar{\bar{\alpha}}<\pi^{\mathcal{M}}$ such that for all $m\in\Box(x*y)$, $(L[r])^{\mathcal{M}}\vDash (S_{\bar{\bar{\alpha}}+m}\;\hbox{is not stationary})$.

\bigskip
Indeed, if  $x,y$ are reals in $L^{\PP_0*\PP^*_{\kappa,\nu,\lambda}}$ such that $x\dot{<}_\pi[G] y$, our bookkeeping guarantees that for some $\xi\in \hbox{Lim}(C^*_0)$, $x=x^\zeta_{{\xi}_0}$, $y=x^\zeta_{\xi_1}$, where $\xi=l_\zeta+\eta$ for some $\eta<\lambda$, and so the real $u_\xi$ codes a
stationarity pattern for $\Box(x*y)$ at $\xi$.

Now, suppose $x,y$ are reals in $L^{\PP_0*\PP^*_{\kappa,\nu,\pi}}$, with the property that for some real $r$, for every countable suitable model $\mathcal{M}$ such that $r\in\mathcal{M}$, there is $\bar{\bar{\alpha}}<\pi^{\mathcal{M}}$ such that for all $m\in\Box(x*y)$, $(L[r])^{\mathcal{M}}\vDash (S_{\bar{\bar{\alpha}}+m}\;\hbox{is not stationary})$. By L\"owenheim-Skolem Theorem, this property holds for arbitrarily large models $\M$
containing the real $r$
and so in particular it holds in $\mathbb{H}_\Theta^{\PP_0*\PP^*_{\kappa,\nu,\lambda}}$, where $\Theta$ is sufficiently large. Thus, there is some $\xi<\pi$,
such that for every $m\in\Box(x*y)$, $L_\Theta[r]\vDash (S_{\xi+m}$ is non-stationary). Since there is no accidental coding of a kill of stationarity by a real,
the ordinal $\xi$ must be in $\hbox{Lim}(C^*_0)$ and so the $u_\xi$ adjoined by $\QQ^*_{\Delta^*(\xi),\xi}$ codes a stationarity pattern for $\Box(a*b)$, where $a <_\pi b$ are the reals given by the bookkeeping function at stage $\xi$ of the iteration. But then $\Box(x*y)\subseteq\Box(a*b)$, which implies that $x=a$, $y=b$ and so indeed $x<_\pi y$.

\section{Discussion and questions}\label{SecQ}

Though the 3D-coherent systems we constructed yield models of several
values in Cicho\'n's diagram, it is still restricted (as in
\cite{mejia}) to constellations where the right side of the diagram
assumes at most 3 different values. So far, the only known model of
more than 3 values on the right (actually 5) is constructed in
\cite{FGKS} with a proper $\omega^\omega$-bounding forcing by a large
product of creatures (though it is restricted to
$\cov(\Nwf)=\dfrak=\aleph_1$).

As discussed before Corollary \ref{Nonewreals}, in all our
constructions we only add two types of generic reals: full generic
reals and restricted generic reals. Different types of generic reals
could be considered (like a real which is restricted generic in some
plane but full generic in the perpendicular plane), but the known
attempts so far destroy the complete embedability of the posets in
the system and, therefore, the construction collapses. Success in this
problem of using a different type of generic real in 3D-coherent
systems would lead to models where more than 3 different values can be
obtained in the right side of Cicho\'n's diagram. For instance,

\begin{question}[{\cite[Sect. 7]{mejia}}]
   Is it consistent with ZFC that $\cov(\Mwf)<\dfrak<\non(\Nwf)<\cof(\Nwf)$?
\end{question}

It seems natural to expect that similar 3D-systems of iterations can be helpful in providing models in
which, for example, $\mathfrak{b}$, $\mathfrak{s}$ and $\mathfrak{a}$ are pairwise distinct. There are three ZFC
admissible constellations: $\sfrak<\bfrak<\afrak$, $\bfrak<\sfrak<\afrak$ and $\bfrak<\afrak<\sfrak$. The consistency
of $\sfrak=\aleph_1<\bfrak<\afrak$ holds in Shelah's original template model~\cite{Sh:700}, while the consistency of
$\aleph_1<\sfrak<\bfrak<\afrak$ has been obtained by the first and third author of the current paper using the iteration
of non-definable (i.e. not Suslin) posets along a Shelah template (see~\cite{VFDM16}). The consistency of $\bfrak<\sfrak<\afrak$
(assuming the existence of a supercompact cardinal) is due to D. Raghavan and S. Shelah~\cite{DRSS2017}, and has been recently
announced at the Oberwolfach Set Theory Meeting, February 2017. Thus, one of the most prominent remaining open questions is the following:

\begin{question}[{\cite[\S 6]{VFJB11}}]
   Is it consistent with ZFC (even assuming large cardinals) that $\bfrak<\afrak<\sfrak$?
\end{question}

We should point out though, that if we are to construct a 3D-system for the above constellation, in order to increase $\sfrak$, we have to include in the construction, non-definable, ccc posets which adjoin non-restricted, unsplitting reals (e.g. we could adjoin full Mathias-Prikry generics). This leads however
to many serious technical problems.




{\small
\bibliography{applmarch2017}
\bibliographystyle{alpha}
}

\end{document}